\documentclass[english]{article}
\usepackage[T1]{fontenc}
\usepackage[latin9]{inputenc}
\usepackage{float}
\usepackage{amstext}
\usepackage{amssymb}
\usepackage{stmaryrd}
\usepackage{graphicx}
\usepackage{esint}

\makeatletter
\newcommand{\lyxaddress}[1]{
	\par {\raggedright #1
	\vspace{1.4em}
	\noindent\par}
}

\@ifundefined{date}{}{\date{}}
\@ifundefined{showcaptionsetup}{}{%
 \PassOptionsToPackage{caption=false}{subfig}}
\usepackage{subfig}
\makeatother

\usepackage{babel}
\begin{document}

\title{Conic Representations of Topological Groups}

\author{Matan Tal}
\maketitle

\lyxaddress{\begin{center}
The Hebrew University of Jerusalem
\par\end{center}}
\begin{abstract}
We define basic notions in the category of conic representations of
a topological group and prove elementary facts about them. We show
that a conic representation determines an ordinary dynamical system
of the group together with a multiplier, establishing facts and formulae
connecting the two categories. The topic is also closely related to
the affine representations of the group. The central goal was attaining
a better understanding of irreducible conic representations of a group,
and - particularly - to determine whether there is a phenomenon analogous
to the existence of a universal irreducible affine representation
of a group in our category (the general answer is negative). Then
we inspect embeddings of irreducible conic representations of semi-simple
Lie groups in some ``regular'' conic representation they possess.
We conclude with what is known to us about the irreducible conic representations
of $SL_{2}\left(\mathbb{R}\right)$.\\
\end{abstract}
\tableofcontents{}

\section{Introduction}

Given a topological group $G$, a finite regular borel measure $\mu$
on $G$ and $\lambda>0$, one may be interested in the measures $\nu$
on $G$ satisfying $\mu*\nu=\lambda\nu$. The set of such $\nu$'s
forms a translation-invariant cone in the linear space of measures
on $G$ (or functions when the $\nu$'s are absolutely continuous
with respect to the Haar measure of the group). In \cite{key-3} Choquet
and Deny assume $G$ is locally compact and commutative and then characterize
the extreme rays of these cones. Furstenberg obtains in \cite{key-2}
an analogous result when $G$ is taken to be a semi-simple Lie group.
It is also discussed in \cite{key-2} which cones of functions yield
irreducible conic representations of the group (in a sense to be made
precise).\\

The above mentioned papers study what may be called the ``regular''
conic representation of the group. Namely conic representations whose
elements are measures or functions on the group, and the group acts
on them by right-translation. Our aim in the present work is to improve
our understanding of conic representations of a topological group
from an abstract point of view.\\

This research was carried out by the author for his master's thesis
done under the supervision of Professor Hillel Furstenberg, who also
suggested the topic. At this point the author wishes to thank him
for his support, patience and generosity.

\section{Affine Representations}

The main source of inspiration of the results to be presented concerning
conic representations was in the theory of affine representations.
We shall therefore begin with summarizing some fundamental ideas of
the latter. The interested reader is referred to {[}2{]} for more
details.\\

The field of \emph{Linear Representations} treats the linear actions
of groups on linear spaces. In analogy to that, \emph{Affine Representations}\textit{\emph{
(the term }}\emph{Affine Dynamics}\textit{\emph{ can be used interchangeably)
}}deals with continuous affine actions of topological groups on compact
convex sets (CCS) in topological linear spaces \footnote{Topological linear spaces are assumed to be Hausdorff throughout this
text.} with a point-separating continuous dual. More precisely, given a
topological group $G$, an \emph{affine representation} of $G$ is
a CCS: $Q\subseteq X$ where $X$ is a topological linear space (either
real or complex) with a point-separating continuous dual\footnote{It is equivalent to the more traditionally written assumption that
$Q$ lies in the dual of some Banach space with the weak-{*} topology,
since it can be naturally embedded in the dual space of the space
of continuous affine functions on $Q.$ It can be done because the
continuous linear functionals separate points.} together with a continuous mapping $\rho:G\times Q\rightarrow Q$,
such that it is an action (i.e. $\rho_{e}=id_{V}$ and $\rho_{g_{1}g_{2}}=\rho_{g_{1}}\circ\rho_{g_{2}}$
for all $g_{1},g_{2}\in G$ ) and $\rho_{g}$ preserves convex combinations
for every fixed $g\in G$. Here and throughout this paper when there
is no danger of misunderstanding, we sometimes refer to such a representation
simply as $Q$ (without mentioning the action explicitly). \\

A morphism in the category of affine systems of a group $G$ is defined
naturally. Given two affine systems $\left(G,Q_{1},\rho_{1}\right)$
and $\left(G,Q_{2},\rho_{2}\right)$, $\varphi:Q_{1}\rightarrow Q_{2}$
is a \emph{morphism} - we call it an \emph{affine homomorphism} -
if it is continuous, affine and commutes with the actions, i.e. $\varphi\circ\left(\rho_{1}\right)_{g}=\left(\rho_{2}\right)_{g}\circ\varphi$
for all $g\in G$. If the morphism is onto, $Q_{2}$ will be called
a \emph{factor} of $Q_{1}$ (this is the quotient mapping of the category).\\

\subsection{Relation to Topological Dynamics}

Given a topological group $G$, a\emph{ topological dynamical system}
is a triple $\left(G,X,\tau\right)$ where $X$ a compact Hausdorff
space, and $\tau:G\times X\rightarrow X$ is continuous. $\left(G,X,\tau\right)$
induces an affine representation $\left(G,\Pr\left(X\right)\right)$,
where $\Pr\left(X\right)$ denotes the space of regular borel probability
measures, and the topology is the weak-{*} topology under the identification
of the measures with bounded linear functionals of $C\left(X\right)$
(either real or complex continuous functions, it does not matter).
Compactness follows from Banach-Alaoglu theorem, which states that
in the dual of a banach space, the closed unit ball (defined by the
norm of the linear functionals) is compact in the weak-{*} topology
(a proof can be found in \cite{key-1}). The action is the push-forward
of the measures: $g_{*}\nu\left(A\right)=\nu\left(g^{-1}\left(A\right)\right)$
for all $g\in G$ , $\nu\in\Pr\left(X\right)$ and Borel sets $A$
of $X$ (if you wish to regard $\nu$ as a linear functional then
every $f\in C\left(X\right)$ satisfies $g\nu\left(f\right)=\nu\left(f\circ g\right)$).
It is affine and continuous. \\

A \emph{homomorphism} between two topological dynamical systems, $\left(G,X_{1},\tau_{1}\right)$
and $\left(G,X_{2},\tau_{2}\right)$, is a continuous mapping $\varphi:X_{1}\rightarrow X_{2}$
which commutes with the actions, i.e. every $g\in G$ satisfies $\varphi\circ\tau_{1}(g)=\tau_{2}(g)\circ\varphi$.
It induces in a natural way an affine homomorphism of affine representations
$\Pr_{*}\varphi=\varphi_{*}:\Pr\left(X_{1}\right)\rightarrow\Pr\left(X_{2}\right)$.
In fact, $\Pr_{*}$ is a covariant functor between the category of
topological dynamical systems of G to the category of affine representations
of G. It is defined either by the duality functor composed on the
pullback of continuous functions, or equivalently, by taking as the
measure of every Borel set of $X_{2}$ the measure of its preimage.
The affine dynamical systems category is itself also a subcategory
of the topological dynamical systems category, and the functor will
also be regarded as a functor from it to itself at will. \\

Conversely, given an affine system $\left(G,Q\right)$, it induces
a topological dynamical system by taking the closure\footnote{The set of extreme points of a CCS is not necessarily closed. A classic
example being the following CCS in $R^{3}$: $Conv\left(\left\{ \left(1,0,\pm1\right)\right\} \cup\left\{ \left(\cos\theta,\sin\theta,0\right):0\leq\theta\leq2\pi\right\} \right)$.} of its extreme points (the extreme points form a $G-invariant$ set).
However, this is not a functor, since an affine homomorphism between
the affine systems $\left(G,Q_{1}\right)$ and $\left(G,Q_{2}\right)$
does not always restrict to a homomorphism between the topological
dynamical systems $\left(G,\overline{ExtQ_{1}}\right)$ and $\left(G,\overline{ExtQ_{2}}\right)$.
As an example, consider the projection of the closed unit disc on
an interval with any group acting trivially (as the identity).

\paragraph*{Reminder: \textmd{Given a compact Hausdorff space $X$, the extreme
points of $\Pr\left(X\right)$ are exactly $\left\{ \delta_{x}:x\in X\right\} $
where $\delta_{x}$ is the Dirac measure of $x\in X$. The mapping
$\delta:X\rightarrow\Pr\left(X\right)$ is continuous, one-to-one
and closed (continuous between compact and Hausdorff spaces) and thus
also a topological embedding.}\protect \\
\textmd{}\protect \\
}

Recalling the fact above, if $X$ is a compact Hausdorff space, then
$\overline{Ext}\left(\Pr\left(X\right)\right)$ is again isomorphic
to it. However, for a CCS $Q$, $\Pr\left(\overline{ExtQ}\right)$
need not be isomorphic to $Q$ . To see this, take $Q$ to be the
unit disc. $\overline{ExtQ}$ is the unit circle and the circle's
$\Pr$ is a CCS, but this time there exists points that are not a
convex combination of two extreme points while in the unit disc every
point is a convex combination of extreme points. \\

\subsection{The Barycenter of a Measure on a CCS}

This is where the point-separating property of the dual to the ambient
space comes into play.

\paragraph*{Reminder (Krein-Milman) \textmd{(a proof can be found in }\cite{key-1}):
\textmd{Given a CCS $Q$ in a topological linear space with a point-separating
continuous dual, then $Q=\overline{Conv\left(ExtQ\right)}$ ($Conv$
is the convex hull of a set).}\protect \\
\textmd{}\protect \\
}

An important concept (and tool) of our subject is the \emph{barycenter}
(center of mass) of a regular borel probability measure defined on
a CCS. With the aid of Krein-Milman theorem on $\Pr\left(Q\right)$
it can be easily verified that the definition of the barycenter is
invariant under isomorphism of affine systems (and as so it is independent
of the embedding).\\

Let $X$ be a topological linear space with a point-separating continuous
dual and $Q\subseteq X$ be a CCS, the \emph{barycentric} mapping$\beta:\Pr\left(Q\right)\rightarrow Q$
is thus defined: given $\nu\in\Pr\left(Q\right)$, $\beta\left(\nu\right)$
is the unique $x\in Q$ satisfying for any $f\in X^{*}$
\[
f\left(x\right)=\int_{Q}f\left(y\right)\,d\nu\left(y\right).
\]
We will write in short $\beta\left(\nu\right)=\int_{Q}y\,d\nu\left(y\right)$
where the meaning is the one above\footnote{This is the Pettis integral of the identity function on $Q$.}.
The uniqueness follows from the point separating property of the continuous
linear functionals. Existence can be proved by deploying the Krein-Milman
theorem on $\Pr\left(Q\right)$. Given $\nu\in\Pr\left(Q\right)$
by the Krein-Milman theorem there exists a net of convex combinations
of Dirac measures on $Q$ that converges to $\nu$, by passing to
a subnet if necessary we may assume that the net of the corresponding
convex combinations of points converges to some $x\in Q$. Since for
any convex combination of Dirac measures the corresponding convex
combination of points in $Q$ is its barycenter, by the definition
of the weak-{*} topology on $\Pr\left(Q\right)$ and the continuity
of any $f\in X^{*}$ (by definition) we deduce that $x$ is the barycenter
of $\nu$.\\

By the proof of existence it is also clear that the equality holds
for any continuous affine $f$ and not only for continuous linear
functionals, thus we deduce that the barycenter does not depend on
the embedding of $Q$ in the linear space.\\

If $\left(G,Q\right)$ is an affine system then $\beta$ is an affine
homomorphism since it is continuous, preserves convex combinations,
and $g\beta\left(\nu\right)=\beta\left(g\nu\right)$ for all $g\in G,\,\nu\in\Pr\left(Q\right)$.
It is also onto (consider Dirac measures), so $Q$ is a factor of
$\Pr\left(Q\right).$\footnote{$\beta$ is in fact a natural transformation between the functor $\Pr$
from the category of affine representations to itself and the category's
identity functor.} Actually, the Krein-Milman theorem implies that the restriction of
$\beta$ to $\Pr\left(\overline{Ext\left(Q\right)}\right)$ is already
onto $Q$.\\

\textbf{Proposition 1:} \textbf{(i)} If $\pi:Q\rightarrow Q'$ is
an affine map between two CCS's which is onto, and $y\in ExtQ'$ then
$Ext\left(\pi^{-1}\left(y\right)\right)\subseteq ExtQ$.\\

\textbf{(ii)} Let $Q$ be a CCS and $z\in Q$ , then $z$ is an extreme
point of $Q$ if and only if for any $\nu\in\Pr\left(Q\right),\,\beta\left(\nu\right)=z$
implies $\nu=\delta_{z}.$\\

\textbf{Proof:} \textbf{(i)} $\pi^{-1}\left(y\right)\neq\emptyset$
because $\pi$ is onto. Let $x\in Ext\left(\pi^{-1}\left(y\right)\right)$.
If $x\notin Ext\left(Q\right)$ then there exists $x_{1},x_{2}\neq x$
in $Q$ such that $x=px_{1}+qx_{2}$ for some $p,q>0,$ $p+q=1$.
$x\in Ext\left(\pi^{-1}\left(y\right)\right)$, so one of them must
not be in $\pi^{-1}\left(y\right)$, assuming without loss of generality
it is $x_{1}$, then $\pi\left(x_{1}\right)\neq y$. However, $p\pi\left(x_{1}\right)+q\pi\left(x_{2}\right)=y$
in contradiction to $y\in ExtQ'$.\\

\textbf{(ii)} The ``if'' part is obvious by taking the points in
the convex combination to be Dirac measures. For the ``only if''
part, let $z\in ExtQ$ and $\nu\in\Pr\left(Q\right)$with $\beta\left(\nu\right)=z$.
By (i), $Ext\left(\beta^{-1}\left(z\right)\right)\subseteq ExtQ$,
so $Ext\left(\beta^{-1}\left(z\right)\right)$ is composed of dirac
measures, and by the definition of $\beta$ it must be equal to $\left\{ \delta_{z}\right\} $.
By Krein-Milman, $\beta^{-1}\left(z\right)=\overline{Conv\left(\left\{ \delta_{z}\right\} \right)}=\delta_{z}$.
~$\blacksquare$\\

\subsection{Irreducible Affine Systems}

In the category of topological dynamical systems, a system that does
not contain any non-empty subsystem (i.e. a closed invariant subset)
besides itself, is called \emph{minimal}. Zorn's Lemma combined with
the ``Finite Intersection Condition'' characterization of compactness
imply that any topological dynamical system contains a minimal subsystem.\\

In analogy to that, in the affine dynamical systems category, we have
the concept of an \emph{irreducible} system, meaning the system has
no non-empty subsystem (i.e. closed convex invariant subset) besides
itself, and it can be deduced, in a similiar way, that any affine
system contains such a subsystem. \\

\textbf{Proposition 2:} An affine system $\left(G,Q\right)$ is irreducible
if and only if all $x\in Q$ satsify $\overline{Gx}\supseteq\overline{ExtQ}$.
\\

\textbf{Proof:} If the system is not irreducible then for some $x\in Q$
the non-empty proper subsystem $\overline{Gx}$ does not contain all
extreme points, since then by Krein-Milman Theorem it will be equal
to $Q$.\\

For the ``only if'' part, we need to prove $\overline{Gx}\supseteq ExtQ$
.

$\Pr\left(\overline{Gx}\right)$ is a subsytem of $\Pr\left(Q\right)$
and as such it is being tranformed by $\beta$ to a non-empty subsytem
of $Q$ , and because $Q$ is irreducible, we get $\beta\left(\Pr\left(\overline{Gx}\right)\right)=Q$.
Thus, for any $z\in ExtQ$ there exists $\nu\in\Pr\left(\overline{Gx}\right)$
such that $\beta\left(\nu\right)=z$, and by prop. 1, $\nu=\delta_{z}$,
so $z\in\overline{Gx}$.~$\blacksquare$\\

\subsection{Strong Proximality and More on Irreduciblity}

Proximality and strong proximality are notions of topological dynamical
systems. We say that a topological dynamical system $\left(G,X\right)$
is \emph{proximal} if for any pair of points $x_{1},x_{2}\in X$ there
exists a net $\left\{ g_{\alpha}\right\} _{\alpha}$ in $G$ such
that both nets $\left\{ g_{\alpha}x_{1}\right\} _{\alpha}$and $\left\{ g_{\alpha}x_{2}\right\} _{\alpha}$converge
to the same point in $X$. Equivalently, $\left(G,X\right)$ is proximal
if and only if any pair of points $x_{1},x_{2}\in X$ satisfies $\overline{G\left(x_{1},x_{2}\right)}\cap\triangle\left(X\times X\right)\neq\emptyset$,
the latter referring to the diagonal $\left\{ \left(x,x\right):x\in X\right\} \subseteq X\times X$.
\\

We say $\left(G,X\right)$ is \emph{strongly proximal} if $\delta_{x}\in\overline{G\nu}$
for all $\nu\in\Pr\left(X\right),\,x\in X$ (or equivalenty, for any
$\nu\in\Pr\left(X\right)$, open set $U\subseteq X$ and $\epsilon>0$
there exists $g\in G$ such that $g_{*}\nu\left(U\right)>1-\epsilon$).
It can be thought of as a notion of ``uniform proximality''. \\

The following propositon shows in particular that strong proximality
implies proximality. It is in in fact somewhat stronger than proximality,
but we will not present here the particular example which proves it
(see \cite{key-4}). \\

\textbf{Proposition 3:} If $\left(G,X\right)$ is a strongly proximal
system then $\overline{G\left(x_{1},...,x_{n}\right)}\supseteq\triangle\left(X^{n}\right)$
for any $x_{1},...,x_{n}\in X$ .\\

\textbf{Proof:} Let $x\in X$ and take $\nu=\frac{\sum_{1}^{n}\delta_{x_{i}}}{n}$.
There exists a net $\left\{ g_{\alpha}\right\} _{\alpha}$in $G$
such that the net $\left\{ g_{\alpha*}\nu\right\} _{\alpha}$converges
to $\delta_{x}$. We will show that all nets $\left\{ g_{\alpha}x_{i}\right\} _{\alpha}$converge
to $x$. Assuming the contrary, w.l.o.g. $\left\{ g_{\alpha}x_{1}\right\} _{\alpha}$does
not converge to $x$. Hence, there exists an open neighborhood $U$
of $x$, and a subset $\left\{ \gamma_{\alpha}\right\} _{\alpha}$of
the index set with $\alpha\leq\gamma_{\alpha}$, such that $\left\{ g_{\gamma_{\alpha}}x_{1}\right\} _{\alpha}\notin U$.\\
$X$ is compact and Hausdorff and hence is normal, so there is an
open neighborhood $V$ of $x$ that satisfies $\overline{V}\subseteq U$,
and by Urysohn's Lemma there exists a continuous function $f:X\rightarrow\left[0,1\right]$
that vanishes on $X\setminus U$ and its value is identically 1 on
$\overline{V}$. But then,\\
$g_{\gamma_{\alpha}}\nu\left(f\right)=\frac{\sum_{1}^{n}\delta_{g_{\gamma_{\alpha}}x_{i}}\left(f\right)}{n}=\frac{\sum_{1}^{n}f\left(g_{\gamma_{\alpha}}x_{i}\right)}{n}=\frac{\sum_{2}^{n}f\left(g_{\gamma_{\alpha}}x_{i}\right)}{n}\leq\frac{n-1}{n}<1=\delta_{x}\left(f\right)$,
in contradiction to to the convergence of $\left\{ g_{\alpha}\nu\right\} _{\alpha}$to
$\delta_{x}$.~$\blacksquare$\\

\textbf{Proposition 4:} If $\left(G,Q\right)$ is an irreducible affine
system then $\left(G,\overline{ExtQ}\right)$ is strongly proximal.\\

\textbf{Proof:} Let $\nu\in\Pr\left(\overline{ExtQ}\right)$, $z\in ExtQ$.
$\nu\in\Pr\left(Q\right)$ so we can take its barycenter in $Q$,
and there is a net $g_{\alpha}$ for which $g_{\alpha}\beta\left(\nu\right)\rightarrow z$
(by prop. 2). $g_{\alpha}\beta\left(\nu\right)=\beta\left(g_{\alpha*}\nu\right)$
and therefore $z\in\beta\left(\overline{G_{*}\nu}\right)$ and by
prop. 1 we get $\delta_{z}\in\overline{G_{*}\nu}$. Since $\overline{G_{*}\nu}$
is closed, the same holds for $z\in\overline{ExtQ}$.~$\blacksquare$\\

The ``converse'' claim also holds.\\

\textbf{Proposition 5:} If $\left(G,X\right)$ is a strongly proximal
dynamical system then $\left(G,\Pr\left(X\right)\right)$ is irreducible.\\

\textbf{Proof:} Follows directly from the definition of strong proximality
and prop. 2.~$\blacksquare$\\

Notice that when restricting to the category of irreducible affine
systems of $G$, $\overline{Ext}$ is in fact a covariant functor
to the category of strongly proximal systems of $G$. For if $\varphi:Q_{1}\rightarrow Q_{2}$
is an affine homomorphism of irreducible affine systems, then $\varphi\left(\overline{ExtQ_{1}}\right)$
is strongly proximal and hence minimal, but it contains $\overline{ExtQ_{2}}$
(by prop. 1 (i)) that is itself minimal and therefore $\varphi\left(\overline{ExtQ_{1}}\right)=\overline{ExtQ_{2}}$.\\

\textbf{Proposition 6:} For any topological group $G$, there exists
an irreducible affine system $\left(G,Q_{G}\right)$, that admits
an affine homomorphism onto any other irreducible affine system of
$G$.\\

\textbf{Proof:} By prop. 2 and Krein-Milman theorem irreducible affine
systems of $G$ are bounded in cardinality by some cardinal $\kappa$,
so if we take $Q_{i}$ to be all irreducible systems of $G$ with
elements in $\kappa$ then there are representatives of every isomorphism
type of irreducible systems. The direct product $\prod Q_{i}$ is
an affine system and thus posseses an irreducible subsystem $Q$ that
by definition admits an affine homomorphism to any other irreducible
system. If each $Q_{i}$ lies in a topological linear space $X_{i}$
with a point-separating continuous dual, then $\prod Q_{i}$ lies
in $\prod X_{i}$ which is also a topological linear space with a
point-separating continuous dual. $\blacksquare$\\

$Q_{G}$ will be called a \emph{universal irreducible affine system}
of $G$. It will soon be shown to be unique and the indefinite article
will be replaced by a definite one.\\

\textbf{Proposition 7:} \textbf{(i)} If $Q_{G}$ is a universal irreducible
affine system of $G$ then $\overline{ExtQ_{G}}$ is a universal strongly
proximal system of $G$ (in the same sense). \textbf{(ii)} If $\Pi_{G}$
is a universal strongly proximal system of $G$ then $\Pr\left(\Pi_{G}\right)$
is a universal irreducible system of $G$.\\

\textbf{Proof: (i)} Let $X$ be a strongly proximal $G$-space. Then
by prop. 5, $\Pr\left(X\right)$ is irreducible, and thus there exists
an affine homomorphism $\tau:Q_{G}\rightarrow\Pr\left(X\right)$ which
in its turn induces $\overline{Ext}_{*}\tau:\overline{ExtQ_{G}}\rightarrow X$.
\\

\textbf{(ii) }Let $Q$ be an irreducible affine $G$-space. Then,
by prop. 4, $\overline{ExtQ}$ is strongly proximal, and thus there
exists a homomorphism $\pi:\Pi_{G}\rightarrow\overline{ExtQ}$ which
in its turn induces $\Pr_{*}\pi:\Pr\left(\Pi_{G}\right)\rightarrow\Pr\left(\overline{ExtQ}\right)$.
Composing the embedding of $\Pr\left(\overline{ExtQ}\right)$ in $\Pr\left(Q\right)$
and then the barycentric mapping finishes the proof.~$\blacksquare$
\\

The next proposition is in the spirit of Schur's Lemma (on irreducible
linear representations).\\

\textbf{Proposition 8:} If $\left(G,Q_{1}\right)$, $\left(G,Q_{2}\right)$
are irreducible affine systems and $\varphi_{1},\varphi_{2}:Q_{1}\rightarrow Q_{2}$
are affine homomorphisms, then $\varphi_{1}=\varphi_{2}$.\\

\textbf{Proof:} $\varphi:=\frac{\varphi_{1}+\varphi_{2}}{2}$ is also
an affine homomorphism. If $\varphi\left(x\right)\in ExtQ_{2}$ then
$\varphi_{1}\left(x\right)=\varphi_{2}\left(x\right)$. $Q_{2}$ is
irreducible, so $\varphi$ is onto, and hence $\left\{ x\in Q_{1}:\varphi_{1}\left(x\right)=\varphi_{2}\left(x\right)\right\} \neq\emptyset$.
But this set is a $G$-invariant CCS and thus is equal to $Q_{1}$.~$\blacksquare$\\

\textbf{Remark:} With a slight modification of the argument, prop.
8 is still valid if we substitute the requirement of irreduciblity
of $Q_{1}$ for the requirement of $\overline{ExtQ_{1}}$ to be minimal.
This more general point of view will be the one with an analogous
proof in the Conic Representations category.\\

The last proposition implies that for any topological group there
is a unique universal irreducible affine representation and a unique
universal strongly proximal space.

\subsection{Amenable Groups and Mostow Groups}

A topological group is said to be \emph{amenable} if all of its affine
actions have a fixed point. If we also assume the group to be locally
compact, Hausdorff and second countable this defintion is equivalent
to other definitions of amenability the reader may know. Given an
affine representation $Q$ of a compact group $K$ and choosing $x_{0}\in Q$,
one may easily see that the barycenter of the push-forward of the
Haar measure of $K$ through the orbit mapping of $x_{0}$ is a fixed
point of the action of $K$. So compact groups are amenable. In addition,
by the Markov-Kakutani fixed-point theorem \cite{key-6} it follows
that abelian groups are amenable. It is also not hard to prove amenability
is preserved by abelian (even amenable) group extensions and thus
all solvable groups are amenable.\\

A topological group $G$ is said to be a \emph{Mostow group} if it
contains a closed amenable sub-group $P$, such that $G/P$ is compact
and $\left(G,G/P\right)$ is strongly proximal. The universal strongly
proximal space of such a group $G$ may be shown to be $G/P$, and
thus its universal irreducible affine representation is $\Pr\left(G/P\right)$.
\\

For either $SL_{n}\left(\mathbb{R}\right)$ or $GL_{n}\left(\mathbb{R}\right)$
, the quotient with their upper-triangular matrices closed sub-group
- called \textit{the Flag Space} - is compact and the group action
on it can be proven to be strongly proximal. In addition, the upper-triangular
matrices form a solvable and thus amenable sub-group. So $SL_{n}\left(\mathbb{R}\right)$
and $GL_{n}\left(\mathbb{R}\right)$ are Mostow groups, and thus the
flag space is their universal strongly proximal space, and its regular
probability measures their universal irreducible representation.\\

\section{The Definition of Conic Representations}

A $cone$ $V$ in a real or complex linear space is a subset of the
space closed under $non-negative\,linear\,combinations$ of its elements,
i.e. $ax+by\in V$ for all $x,y\in V$ and $a,b\geq0$ \footnote{Some authors use the term $convex\,cone$ for describing such a cone.}.
The sets $\left\{ \lambda x|\lambda>0\right\} $, where $0\neq x\in V$,
are called $rays$ of the cone. A ray $\left\{ \lambda x|\lambda>0\right\} $
is called an $extreme\,ray$ if $y_{1},y_{2}\in V$, $a,b>0$, $ay_{1}+by_{2}=x$
implies $y_{1},y_{2}\in\left\{ \lambda x|\lambda\geq0\right\} $.
A $conic\,function$ from a cone to $\mathbb{R}$ or $\mathbb{C}$
is a function preserving non-negative linear combinations. \\

Let $V$ be a cone contained in a (real or complex) topological linear
space $X$ with a point-separating continuous dual. A subset $Q\subseteq V$
is called a $section$ of $V$ if there exists a continuous conic
function $L:V\rightarrow\mathbb{R}_{\geq0}$ such that all $x\in V\setminus0$
satisfy $L\left(x\right)>0$ and $Q=L^{-1}\left(1\right)$. Every
section admits a canonical projection onto itself from the cone. Its
restriction to another section defined by a conic function $L_{1}$
is easily seen to be affine if and only if $L_{1}=aL$ for some $a>0$.
If $V_{1}$, $V_{2}$ are two cones with sections $Q_{1}$ defined
by $L_{1}$ and $Q_{2}$ defined by $L_{2}$ respectively, and $\phi:Q_{1}\rightarrow Q_{2}$
is an affine (preserving convex combinations) homeomorphism, then
$\varphi:V_{1}\rightarrow V_{2}$ defined by $\varphi\left(x\right)=L_{1}\left(x\right)\phi\left(\frac{x}{L_{1\left(x\right)}}\right)$
is an isomorphism of the two cones - i.e. $\psi\left(y\right)=L_{2}\left(y\right)\phi^{-1}\left(\frac{y}{L_{2\left(y\right)}}\right)$
is its inverse, and they are both continuous and preserve non-negative
linear combinations. In short, two cones having isomorphic sections
are isomorphic, so we can speak of $the$ cone $V_{Q}$ $generated$
by the section $Q\subseteq X$. Every CCS $Q$ can be viewed as a
section of some cone since the cone in $X\oplus\mathbb{R}$ generated
by $\left\{ \left(x,1\right)|x\in Q\right\} $ is such a cone.\\

For our needs we require $V$ to be non-zero and to admit a compact
section\footnote{If a section is compact, then all sections are, since they are homeomorphic
by continuity of the functions defining them. It is unknown to me
whether all sections of a such a cone are affinely isomorphic, but
the important fact is that the natural projection of one's extreme
points to another's need not have an extension to an affine isomorphism.}. $V$ is thus closed, Hausdorff and locally compact. In this setting
- and given a topological group $G$ - we define a $conic\,representation$
of $G$ on $V$ to be a continuous mapping $\tau_{\left(\cdot\right)}\left(\cdot\right):G\times V\rightarrow V$
such that it is an action (i.e. $\tau_{e}=id_{V}$ and $\tau_{g_{1}g_{2}}=\tau_{g_{1}}\circ\tau_{g_{2}}$
for all $g_{1},g_{2}\in G$ ) and $\tau_{g}$ preserves non-negative
linear combinations for any $g\in G$ . It may also be called a $conic\,dynamical\,system$,
and $V$ may also be called a $G-cone$.\\

\textbf{Example 3.1:} The action of the group of matrices $\left\{ \left(\begin{array}{cc}
a & 0\\
0 & b
\end{array}\right):a,b>0\right\} \cup\left\{ \left(\begin{array}{cc}
0 & a\\
b & 0
\end{array}\right):a,b>0\right\} $ on the closed first quadrant of $\mathbb{R}^{2}$.\\

A $sub-representation$ is a closed (non-zero) sub-cone invariant
under the action of $G$. In the example there are no non-empty sub-representations
properly contained in the original since the action is transitive
on the interior. \\

\section{Conic Pairs}

Given a topological $G-$space $M$ with a continuous action $\rho_{\left(\cdot\right)}\left(\cdot\right):G\times M\rightarrow M$.
A $multiplier$ of $G$ and $M$ (relative to $\rho$) is a continuous
function $\mbox{\ensuremath{\sigma}:}G\times M\rightarrow\mathbb{R}_{>0}$
that satisfies $\sigma\left(g\gamma,x\right)=\sigma\left(g,\rho_{\gamma}\left(x\right)\right)\sigma\left(\gamma,x\right)$
for all $g,\gamma\in G$, $x\in M$ (in particular, $\sigma\left(e,\cdot\right)\equiv1$)
\footnote{As a function from $G$ to the $G$-module of functions on $M$ it
is also referred to as a crossed homomorphism or a 1-cocycle (see
\cite{key-2}).}. Any continuous homomorphism $\varphi:G\rightarrow\mathbb{R}_{>0}$
induces in a natural way a multiplier not depending on $M$.\\

If $\left(G,V,\tau\right)$ is a conic representation and $Q\subseteq V$
is a compact section, then since the action of any $g\in G$ maps
a ray to a ray, it induces a continuous action of $G$ on $Q$. Denoting
it $\rho_{\left(\cdot\right)}\left(\cdot\right):G\times Q\rightarrow Q$,
there exists a unique function $\mbox{\ensuremath{\sigma}:}G\times Q\rightarrow\mathbb{R}_{>0}$
satisfying $\tau_{g}\left(x\right)=\sigma\left(g,x\right)\rho_{g}\left(x\right)$
for all $g\in G$, $x\in Q$ . Taking $L$ to be the continuous function
defining $Q$, as above, we get

\[
\sigma\left(g,x\right)=L\left(\tau_{g}\left(x\right)\right).
\]
\\

In particular, $\sigma$ is continuous. \footnote{This is an instance of a more general phenomenon in dynamics. Given
groups $G$ and $H$, a $G$-space $X$, a $H$-space $M$ and a multiplier
$\sigma:G\times X\rightarrow H$, then $X\times M$ can be taken to
be a $G$-space under the action $g\left(x,m\right)=\left(gx,\sigma\left(g,x\right)m\right)$
for any $g\in G$, $x\in X$, $m\in M$ (the fact that this is an
action is equivalent to $\sigma$ having the multiplier property).
With the action of $G$ thus defined $X\times M$ is called a skew
product with multiplier $\sigma$.} \\

\textbf{Proposition 9:} $\sigma$ is a multiplier of $G$ and $Q$.\\

\textbf{Proof:} Only the last condition in the multiplier definition
requires some calculation: $\tau_{g\gamma}\left(x\right)=\sigma\left(g\gamma,x\right)\rho_{g\gamma}\left(x\right)$
but it also equals $\tau_{g}\circ\tau_{\gamma}\left(x\right)=\tau_{g}\left(\sigma\left(\gamma,x\right)\rho_{\gamma}\left(x\right)\right)=\sigma\left(\gamma,x\right)\tau_{g}\left(\rho_{\gamma}\left(x\right)\right)=\sigma\left(\gamma,x\right)\sigma\left(g,\rho_{\gamma}\left(x\right)\right)\rho_{g}\left(\rho_{\gamma}\left(x\right)\right)=\sigma\left(\gamma,x\right)\sigma\left(g,\rho_{\gamma}\left(x\right)\right)\rho_{g\gamma}\left(x\right)$
and comparing coefficients we are done. $\blacksquare$\\

\textbf{Proposition 10:} $\sigma\left(g,\lambda x_{1}+\left(1-\lambda\right)x_{2}\right)=\lambda\sigma\left(g,x_{1}\right)+\left(1-\lambda\right)\sigma\left(g,x_{2}\right)$
for all $0\le\lambda\le1$, $g\in G$, $x_{1},x_{2}\in Q$. In particular,
this implies that for each $g\in G$, $\sigma\left(g,\cdot\right)$
is determined by its values on $ExtQ$.\\

\textbf{Proof:} $\tau_{g}\left(\lambda x_{1}+\left(1-\lambda\right)x_{2}\right)=\lambda\tau_{g}\left(x_{1}\right)+\left(1-\lambda\right)\tau_{g}\left(x_{2}\right)$\\

$\sigma\left(g,\lambda x_{1}+\left(1-\lambda\right)x_{2}\right)\rho_{g}\left(\lambda x_{1}+\left(1-\lambda\right)x_{2}\right)=\lambda\sigma\left(g,x_{1}\right)\rho_{g}\left(x_{1}\right)+\left(1-\lambda\right)\sigma\left(g,x_{2}\right)\rho_{g}\left(x_{2}\right)$\\

$\rho_{g}\left(\lambda x_{1}+\left(1-\lambda\right)x_{2}\right)=\frac{\lambda\sigma\left(g,x_{1}\right)}{\sigma\left(g,\lambda x_{1}+\left(1-\lambda\right)x_{2}\right)}\rho_{g}\left(x_{1}\right)+\frac{\left(1-\lambda\right)\sigma\left(g,x_{2}\right)}{\sigma\left(g,\lambda x_{1}+\left(1-\lambda\right)x_{2}\right)}\rho_{g}\left(x_{2}\right)$~~~~$\left(1\right)$\\

The coefficients on the right side of the equation need to sum up
to 1 in order for the linear combination to stay in $Q$, and that
finishes the proof. $\blacksquare$\\

\textbf{Remarks:} 
\begin{itemize}
\item The last propostion could have also been proved by using the equality
$\sigma\left(g,x\right)=L\left(\tau_{g}\left(x\right)\right)$, but
we wanted to also deduce $\left(1\right)$.
\item $\left(1\right)$ can be generalized to the statement that $\rho_{g}\left(\beta\left(\nu\right)\right)=\int\frac{\sigma\left(g,x\right)}{\sigma\left(g,\beta\left(\nu\right)\right)}\rho_{g}\left(x\right)d\nu\left(x\right)$
for any $\nu\in\Pr\left(Q\right)$. The integration is Pettis integration,
as in the definition of the barycenter, and the proof follows by considering
convex combinations of Dirac measures and then using the Krein-Milman
theorem.
\item By $\left(1\right)$, $\rho_{g}$ transfers straight lines into straight
lines. In the terminology of projective geometry it is a $collineation$.
This is not surprising because it is defined in an analogous manner
to a projective transformation: given a plane in a linear space $X$
defined as the level set $L\equiv1$ of a linear functional $L$,
a projective transformation is any transformation from the plane to
itself defined by $\frac{T\left(x\right)}{L\left(T\left(x\right)\right)}$
where $T:X\rightarrow X$ is a linear isomorphism. Thus we shall call
an action $\rho$ of $G$ on a CCS $Q$ that can be obtained as the
induced action on a section of some conic representation a $projective\,action$
of $G$ on $Q$\footnote{Note that in the theory of linear representations, a projective representation
of a topological group $G$ in $\mathbb{PR}^{n}$ is a continuous
homomorphism $G\rightarrow PGL_{n}\left(\mathbb{R}\right)$. A stronger
(non-equivalent!) condition is that this homomorphism factors continuously
through $GL_{n}\left(\mathbb{R}\right)\rightarrow PGL_{n}\left(\mathbb{R}\right)$.
The definition we have just given for an action on a CCS is analogous
to the latter.}. 
\end{itemize}
A natural question that arises is about reversing the standpoint of
the previous analysis: starting with a CCS $Q$ in a topological linear
space with a point-separating continuous dual\footnote{From now on assumed without further remark.},
$\rho_{\left(\cdot\right)}\left(\cdot\right):G\times Q\rightarrow Q$
a continuous action of $G$ on $Q$, and a multiplier $\mbox{\ensuremath{\sigma}:}G\times Q\rightarrow\mathbb{R}_{>0}$
(relative to $\rho$), when could they be synthesized into a conic
representation of $G$ on $V_{Q}$ inducing the action $\rho$ on
$Q$?\\

If there exists such a representation $\tau$, it is necessarily unique
because $\tau_{g}(\lambda x)=\lambda\tau_{g}\left(x\right)=\lambda\sigma\left(g,x\right)\rho_{g}\left(x\right)$
for all $g\in G$, $\lambda\geq0$, $x\in Q$. Is $\tau$ (thus defined)
always a conic representation? It is a continuous action. So a necessary
and sufficient condition for $\tau$ to be a conic representation
is that it preserves non-negative linear combinations and this is
equivalent, for a (non-negative) homogenous $\tau$ (like ours), to
preserving convex combinations; hence to $\left(1\right)$ (the three
equations in the proof of prop. 10 are equivalent to each other).
We summarize the conclusions of the last two paragraphs in the following
proposition.\\

\textbf{Proposition 11:} Given a CCS $Q$, $\rho_{\left(\cdot\right)}\left(\cdot\right):G\times Q\rightarrow Q$
a continuous action of $G$ on $Q$ together with a multiplier $\mbox{\ensuremath{\sigma}:}G\times Q\rightarrow\mathbb{R}_{>0}$,
then $\rho$ and $\sigma$ can be induced by a conic representation
of $G$ if and only if they together satisfy $(1)$. In this case,
the conic representation is unique (up to isomorphism).\\

We call a pair $\left(\rho,\sigma\right)$ that satisfies $\left(1\right)$
a $conic\,pair$ of $G$ on $Q$, and $\sigma$ a $conic\,multiplier$
of $\rho$. So an action $\rho$ on a CCS $Q$ is a projective action
if and only if there exists a multiplier $\sigma$ such that together
they form a conic pair.\\

\textbf{Proposition 12:} If $\left(\rho,\sigma\right)$, $\left(\rho,\sigma'\right)$
are conic pairs of $G$ on a CCS $Q$ then:\\

$\left(i\right)$ There exists a unique function $a:G\rightarrow\mathbb{R}_{>0}$
such that $\sigma'\left(g,x\right)=a\left(g\right)\sigma\left(g,x\right)$
for all $g\in G$ , $x\in Q$ .\\

$\left(ii\right)$ $a$ is a continous homomorphism.\\

\textbf{Proof: }

\textbf{$\left(i\right)$}

$\frac{\lambda\sigma\left(g,x_{1}\right)}{\sigma\left(g,\lambda x_{1}+\left(1-\lambda\right)x_{2}\right)}\rho_{g}\left(x_{1}\right)+\frac{\left(1-\lambda\right)\sigma\left(g,x_{2}\right)}{\sigma\left(g,\lambda x_{1}+\left(1-\lambda\right)x_{2}\right)}\rho_{g}\left(x_{2}\right)=\frac{\lambda\sigma'\left(g,x_{1}\right)}{\sigma'\left(g,\lambda x_{1}+\left(1-\lambda\right)x_{2}\right)}\rho_{g}\left(x_{1}\right)+\frac{\left(1-\lambda\right)\sigma'\left(g,x_{2}\right)}{\sigma'\left(g,\lambda x_{1}+\left(1-\lambda\right)x_{2}\right)}\rho_{g}\left(x_{2}\right)$\\

Assuming $x_{1}\neq x_{2}$, then $\rho_{g}\left(x_{1}\right)\neq\rho_{g}\left(x_{2}\right)$
since $\rho_{g}$ is invertible. Comparing coefficients and dividing
the two equations one obtains the equality\\

$\frac{\sigma\left(g,x_{1}\right)}{\sigma\left(g,x_{2}\right)}=\frac{\sigma'\left(g,x_{1}\right)}{\sigma'\left(g,x_{2}\right)}$\\

so we take $a\left(g\right)=\frac{\sigma'\left(g,x_{1}\right)}{\sigma\left(g,x_{1}\right)}$.\\

$\left(ii\right)$\\

$a\left(g\gamma\right)\sigma\left(g\gamma,x\right)=\sigma'\left(g\gamma,x\right)=\sigma'\left(\gamma,x\right)\sigma'\left(g,\rho_{\gamma}\left(x\right)\right)=a\left(\gamma\right)a\left(g\right)\sigma\left(\gamma,x\right)\sigma\left(g,\rho_{\gamma}\left(x\right)\right)=a\left(\gamma\right)a\left(g\right)\sigma\left(g\gamma,x\right)$.
$\blacksquare$\\

\textbf{Proposition 13:} If $\left(\rho,\sigma\right)$ is a conic
pair of $G$ on a CCS $Q$ and $\sigma'\left(g,x\right)=a\left(g\right)\sigma\left(g,x\right)$
where $a:G\rightarrow\mathbb{R}_{>0}$ is a homomorphism, then $\left(\rho,\sigma'\right)$
is also a conic pair of $G$ on $Q$.\\

\textbf{Proof:} $\sigma'$ is a multiplier by the calculation in part
$\left(ii\right)$ of the previous proof. The rest is obvious. $\blacksquare$\\

The last two propositions yield the following corollary:\\

\textbf{Corollary 14:} Given a conic pair $\left(\rho,\sigma\right)$
of $G$ on a CCS $Q$, there is a natural one-to-one correspondence
between the conic multipliers of $\rho$ and the continuous homomorphims
from $G$ to $\mathbb{R}_{>0}$ ($\sigma$ corresponds to the trivial
homomorphism).\footnote{$\rho$ can also have no conic pairs at all.}
$\blacksquare$\\

\textbf{Remark:} Homomorphisms from a group to the multiplicative
group $\mathbb{R}_{>0}$ are in one-to-one corresponce with its homomorphisms
to $\mathbb{R}$ by composing with the logarithm function.\\

Given a CCS $Q$, $\sigma\equiv1$ is a multiplier for any $\rho$
of $G$. It forms a conic pair with $\rho$ if and only if $\rho$
is affine (i.e. preserves convex combinations). We call the conic
representation induced by it and an affine $\rho$ a $degenerate\,conic\,representation$.
Note that a conic representation is degenerate if and only if it admits
an invariant section. If $Q$ is an affine representation we call
the conic representation induced by it and $\sigma\equiv1$ $the\,degenerate\,conic\,representation\,that\,belongs\,to\,Q$.\\

\textbf{Example 4.1:} A continuous action of $G$ on a compact Hausdorff
space $X$, induces a conic action on the cone of finite regular measures
of $X$ (The linear space of finite regular measures is identified
with the dual space of $C\left(X\right)$ equipped with the weak-{*}
topology). This conic representation of $G$ is degenerate since $\Pr\left(X\right)$
is an invariant section.\\

\textbf{Corollary 15:} If $\rho$ is an affine representation of $G$,
then $\left(\rho,\sigma\right)$ is a conic pair if and only if $\sigma\left(g,x\right)=a\left(g\right)$
where $a:G\rightarrow\mathbb{R}_{>0}$ is a continous homomorphism
(for $\sigma\equiv1$ forms a conic pair together with $\rho$). $\blacksquare$\\

\section{Homomorphisms of Conic Representations\protect \\
}

A $homomorphism\,\varphi\,between\,two\,conic\,representations$ of
$G$ is defined to be a continuous \textbf{nowhere-zero} mapping,
preserving non-negative linear combinations (a.k.a. $conic\,mapping$)
that commutes with the group action. If $\varphi:\left(V_{1},\tau\right)\rightarrow\left(V_{2},\eta\right)$
is a homomorphism of conic representations which is \textbf{onto}
- this is the quotient mapping in the category of conic representations
- and we say $\left(V_{2},\eta\right)$ is a $factor$ of $\left(V_{1},\tau\right)$.\\

\textbf{Remark:} If a conic representation has a degenerate factor
than it is itself degenerate (the inverse image of an invariant section
is an invariant section). \\

\textbf{Proposition 16:} Let $\varphi:\left(V_{1},\tau\right)\rightarrow\left(V_{2},\eta\right)$
be a homomorphism of conic representations of $G$, $Q_{2}$ a section
of $V_{2}$ and $Q_{1}:=\varphi^{-1}\left(Q_{2}\right)$. If $\left(\rho,\sigma_{1}\right)$,$\left(\theta,\sigma_{2}\right)$
are the conic pairs of $Q_{1}$ and $Q_{2}$ respectively, then $\sigma_{2}\left(g,\varphi\left(x\right)\right)=\sigma_{1}\left(g,x\right)$
for all $g\in G$, $x\in Q_{1}$. Moreover, the restriction $\tilde{\varphi}=\varphi|_{Q_{1}}:Q_{1}\rightarrow Q_{2}$
satisfies $\tilde{\varphi}\circ\rho_{g}=\theta_{g}\circ\tilde{\varphi}$.\\

\textbf{Proof:} First, if $L_{2}$ defines the section $Q_{2}$ then
$L_{1}=L_{2}\circ\varphi$ defines $Q_{1}$ and so it is a section.\\

Let $g\in G,$$x\in Q_{1}$.\\

$\varphi\circ\tau_{g}\left(x\right)=\varphi\left(\sigma_{1}\left(g,x\right)\rho_{g}\left(x\right)\right)=\sigma_{1}\left(g,x\right)\varphi\left(\rho_{g}\left(x\right)\right)$\\

On the other hand,\\

$\eta_{g}\circ\varphi\left(x\right)=\sigma_{2}\left(g,\varphi\left(x\right)\right)\theta_{g}\left(\varphi\left(x\right)\right)$\\

and so,\\

$\sigma_{1}\left(g,x\right)\varphi\left(\rho_{g}\left(x\right)\right)=\sigma_{2}\left(g,\varphi\left(x\right)\right)\theta_{g}\left(\varphi\left(x\right)\right)$\\

Since $\varphi\left(\rho_{g}\left(x\right)\right),\theta_{g}\left(\varphi\left(x\right)\right)\in Q_{2}$
the coefficients are equal. $\blacksquare$\\
\\

\section{The Resultant of a Compactly Supported Measure on a Cone}

Let $V$ be a cone - in a topological linear space $X$ with a point-separating
continuous dual - admitting compact sections. It is thus locally compact.
Let $M_{C}\left(V\right)$ be the space of non-negative regular measures
on $V$ which are compactly supported. Then the resultant $r:M_{C}\left(V\right)\rightarrow V$
is defined by\\

\[
r\left(\mu\right):=\int_{V}yd\mu\left(y\right)
\]

where, as in the definition of the barycenter, the integration is
Pettis integration. The justification for this definition is similar
to the one given in the definition of the barycenter. If there exists
$r\left(\mu\right)$ in $X$ such that $\varphi\left(r\left(\mu\right)\right)=\int_{V}\varphi\left(y\right)d\mu\left(y\right)$
for all continuous linear functionals $\varphi$, then it is the unique
element of $X$ satisfying this since the continuous linear functionals
of $X$ separate points. To see there exists such a $r\left(\mu\right)$
in $V$, one notices that it exists for positve linear combinations
of Dirac measures and that such measures are dense in $M_{C}\left(V\right)$(by
using the Krein-Milman theorem), and continues as in the proof of
existence of the barycenter. The proof of existence of the resultant
implies the required equality holds not only for continuous linear
functionals, but for all continuous conic functions; thus the definition
is an isomorphism invariant of cones. The space $M_{C}\left(V\right)$
is itself a cone and $r$ preserves non-negative linear combinations.
Equipped with the weak-{*} topology on $M_{C}\left(V\right)$, $r$
is \textbf{not} \textbf{continuous}, and $M_{C}\left(V\right)$ does
not admit a compact section. However, if $V$ is a conic representation
of $G$, then $G$ has also a natural action on $M_{C}\left(V\right)$
and the actions do commute with $r$. We may take subcones of $M_{C}\left(V\right)$
that do admit a compact section and that $r$ restricted to them is
continuous. This will be done in section 8 (Semi-Conic Representations).

\section{Irreducible Conic Representations}

A conic representation is called $irreducible$ if it has no non-empty
sub-representations other than itself. Using Zorn's Lemma and compactness
of the section one can show that any conic representation admits an
irreducible sub-representation. \\

We have already seen an example of an irreducible conic representation
(example 3.1). Another class of (degenerate) examples can be obtained
by taking the action of $G$ on $X$ in example 4.1 to be strongly
proximal. However, example 3.1 teaches us that the induced action
of $G$ on closure of the set of extreme points of a section need
not even be proximal (though it necessarily has to be minimal by lemma
17). As we shall see later (in section 10 about $SL_{2}\left(\mathbb{R}\right)$)
the converse is also false. Namely, if the action on the closure of
the set of extreme points of a section is strongly proximal, it does
not guarantee that the conic representation is irreducible.\\

Clearly, a homomorphism of conic representations whose range is irreducible
is necessarily onto.\\

A projective action $\rho$ of $G$ on $Q$ is called $irreducible$
if it has no non-trivial closed convex invariant subset of $Q$ (the
definition is independent of $\sigma)$. Given a conic representation
$V$ and one of its sections $Q$ with an induced action $\rho$,
$V$ is irreducible if and only if $\rho$ is irreducible.\\

\textbf{Lemma 17:} $\rho$ is an irreducible projective action of
$G$ on $Q$ if and only if $\overline{\rho_{G}x}:=\overline{\left\{ \rho_{g}x:g\in G\right\} }$
contains $\overline{Ext\left(Q\right)}$ for all $x\in Q$.\\

\textbf{Proof:} The ``if'' part is trivial. For the ``only if''
part consider $Conv\left(\overline{\rho_{G}x}\right)$. It is an invariant
set under the action $\rho$, since - as already mentioned -$\left(1\right)$
(in the proof of prop. 10) is equivalent to the statement $\rho_{g}\left(\beta\left(\nu\right)\right)=\int\frac{\sigma\left(g,x\right)}{\sigma\left(g,\beta\left(\nu\right)\right)}\rho_{g}\left(x\right)d\nu\left(x\right)$
for any $\nu\in\Pr\left(Q\right)$. Hence given a convex combination
$\lambda_{1}x_{1}+...+\lambda_{n}x_{n}$ where $x_{1},...,x_{n}\in\overline{\rho_{G}x}$,
taking $\nu=\lambda_{1}\delta_{x_{1}}+...+\lambda_{n}\delta_{x_{n}}$we
deduce that $\rho_{g}\left(\lambda_{1}x_{1}+...+\lambda_{n}x_{n}\right)$
is a convex combination of $\rho_{g}\left(x_{1}\right),...,\rho_{g}\left(x_{n}\right)\in\overline{\rho_{G}x}$.\\

Since $Conv\left(\overline{\rho_{G}x}\right)$ is an invariant set,
so is $\beta\left(\Pr\left(\overline{\rho_{G}x}\right)\right)=\overline{Conv\left(\overline{\rho_{G}x}\right)}$,
and since $\rho$ is an irreducible action of $G$ on $Q$ we have
$\beta\left(\Pr\left(\overline{\rho_{G}x}\right)\right)=Q$. By prop.
1, $\delta_{z}\in\Pr\left(\overline{\rho_{G}x}\right)$ for any $z\in ExtQ$,
and thus $z\in\overline{\rho_{G}x}$. $\blacksquare$\\

\textbf{Proposition 18:} If $\varphi_{1},\varphi_{2}:V_{1}\rightarrow V_{2}$
are homomorphisms of conic representations of $G$, with $V_{2}$
irreducible, and the induced action of $G$ on the closure of the
extreme points set of a section $Q$ (and hence all sections) of $V_{1}$
is minimal, then $\varphi_{2}=a\varphi_{1}$ for some $a>0$.\\

\textbf{Proof:} $\varphi:=\varphi_{1}+\varphi_{2}$ is also a homomorphism.
Let $Q_{2}$ be a section of $V_{2}$ and $Q_{1}=\varphi^{-1}\left(Q_{2}\right)$.
If $y_{0}$ belongs to an extreme ray of $V_{2}$, then - since $\varphi$
is onto ($V_{2}$ is irreducible) - there exists $x_{0}\in\varphi^{-1}\left(y_{0}\right)$
that belongs to an extreme ray of $V_{1}$ ($\varphi|_{Q_{1}}:Q_{1}\rightarrow Q_{2}$
is an affine mapping which is onto - now see prop. 1 (i)). Thus $\varphi_{1}\left(x_{0}\right)$
and $\varphi_{2}\left(x_{0}\right)$ also belong to the extreme ray
of $y_{0}$ and $\varphi_{2}\left(x_{0}\right)=a\varphi_{1}\left(x_{0}\right)$
for some $a>0$. This implies that $\varphi_{2}\left(x\right)=a\varphi_{1}\left(x\right)$
for all $x$ in the orbit of $x_{0}$. Since the induced action of
$\overline{ExtQ_{1}}$ is minimal we have $\varphi_{2}\left(x\right)=a\varphi_{1}\left(x\right)$
for all $x\in\overline{ExtQ_{1}}$, but $\varphi_{1},\varphi_{2}$
are affine and the equality holds for all $x\in Q_{1}$. $\blacksquare$\\

\textbf{Corollary 19:} If $\varphi:V_{1}\rightarrow V_{2}$ is a homomorphism
between irreducible conic representations of $G$, then it is unique
up to a multiplication by a positive scalar.\\

In the category of affine representations we have a notion of the
universal irreducible affine representation of a group that always
exists. It means that any other irreducible affine representation
is its factor via a unique homomorphism. \textbf{Does a universal
irreducible conic representation exist for every topological group}
$G$\textbf{ }(i.e. an irreducible representation admitting a homomorphism
to every other one)\textbf{?} Via Corollary 19, we know that if it
exists it is essentially unique.\\

The degenerate conic representation of an irreducible affine representation
of $G$ is an irreducible conic representation, and hence if there
exists a universal irreducible conic representation of $G$ it is
degenerate (consider the homomorphism from it to any degenerate irreducible
representation and the remark before prop. 16). Of course, if there
exist only degenerate irreducible conic representations then the degenerate
one belonging to the universal irreducible affine representation is
the universal irreducible conic representation. This points to a more
general fact.\\

\textbf{Theorem 20:} If $G$ admits a universal irreducible conic
representation $\left(V_{G},\tau\right)$ then it is the degenerate
one belonging to the universal irreducible affine representation.\\

\textbf{Proof:} We have just explained why it is degenerate. If $Q_{G}$
is the universal irreducible affine, let us denote by $\left(\tilde{V},\eta\right)$
the conic representation generated by it. From universality of $V_{G}$,
there exists a homomorphism $\varphi:V_{G}\rightarrow\tilde{V}$.
$Q:=\varphi^{-1}\left(Q_{G}\right)$ is a section of $V_{G}$ with
$\sigma\equiv1$, so the mapping $\varphi|_{Q}$ between the invariant
sections $Q$ and $Q_{G}$ (with the induced actions on them being
affine) is an affine homomorphism of irreducible affine representations,
and from universality of $Q_{G}$ it is an isomorphism. $\blacksquare$\\

The next proposition puts things a bit more in place.\\

\textbf{Propositon 21:} If $Q_{1}$ and $Q_{2}$ are both invariant
sections of a degenerate irreducible conic representation $V$ of
$G$ (with action $\tau:G\times V\rightarrow V$), then one is a multiplication
by a positive scalar of the other. In particular, this implies that
they are isomorphic as affine representations of $G$.\\

\textbf{Proof:} Say $Q_{1}$ is defined in $V$ by $L$. Taking some
$x_{0}\in Q_{2}$, $L\left(\frac{\tau_{g}\left(x_{0}\right)}{L\left(x_{0}\right)}\right)=1$
for all $g\in G$, and thus the set $\overline{\tau_{G}\left(x_{0}\right)}=\overline{\left\{ \tau_{g}\left(x_{0}\right):g\in G\right\} }\subseteq Q_{2}$
has constant $L$ value $a>0.$ Since $V$ is irreducible, then so
is $Q_{2}$, and so $\overline{\tau_{G}\left(x_{0}\right)}$ contains
all extreme points of $Q_{2}$ (lemma 17), and therefore, by Krein-Milman
theorerm, $L\left(x\right)=a$ for all $x\in Q_{2}$. $\blacksquare$\\

\textbf{Lemma 22:} If $\left(V',\eta\right)$ is a degenerate conic
representation, $\varphi:V'\rightarrow V$ is a homomorphism of conic
representations \textbf{onto} $\left(V,\tau\right)$, and $Q$ is
a section of $V$ (defined by $L\equiv1$) with multiplier $\sigma$,
then there exist $0<a<b$ such that $a<\sigma\left(g,x\right)<b$
for all $g\in G$, $x\in Q$.\\
\textbf{Remark:} In the next section we present an example (example
8.1) in which a factor of a degenerate conic representation is not
degenerate. However, we prove that if it is irreducible then it has
to be degenerate (Theorem 34).\\

\textbf{Proof:} Let $Q'$ be a section of $V'$ with multiplier $\sigma'\equiv1$.
Since $\varphi$ is onto, $\varphi\left(Q'\right)$ intersects all
rays of $V$, and it is invariant under the action of $G$ on $V$
since $Q'$ is an invariant section of $V'$. $\varphi\left(Q'\right)$
is compact, and thus $L\left(\varphi\left(Q'\right)\right)\subseteq\left(s,m\right)$
for some $0<s<m$. Thus, $\frac{s}{m}<L\left(\tau_{g}\left(x\right)\right)=\sigma\left(g,x\right)<\frac{m}{s}$
for any $x\in Q$. $\blacksquare$\\

\textbf{Corollary 23:} Let $a:G\rightarrow\mathbb{R}_{>0}$ be a non-trivial
continuous group homomorphism. Let $\left(\rho,\sigma\right)$ be
a conic pair on a CCS $Q$, where $\rho$ is affine and $\sigma\left(g,x\right)=a\left(g\right)$,
and let $\left(V,\tau\right)$ be its generated conic representation.
Then $\left(V,\tau\right)$ is not the range of any conic representation
homomorphism whose domain is degenerate.\\

\textbf{Corollary 24:} If $G$ admits a non-trivial continuous group
homomorphism to $\mathbb{R}_{>0}$ then it has no universal irreducible
conic representation.\\

In particular, if $G$ is locally compact (and Hausdorff) but not
unimodular, the modular function is such a homomophism and thus $G$
has no universal irreducible conic representation.\\

A topological group is $Tychonoff$ if all its conic representations
have invariant rays. It obviously implies amenability, and is in fact
stronger. The group of transformations of the plane that is generated
by rotations and translations is an example of a group which is amenable
but not Tychonoff (see \cite{key-2}).\\

The irreducible conic representations of a Tychonoff group $G$ are
only its one dimensional ones, i.e. its conic actions on the cone
$\mathbb{R}_{\geq0}$. So it admits a universal irreducible conic
representation if and only if its only irreducible conic representation
is the identity action on $\mathbb{R}_{\geq0}$. But the conic actions
of $G$ on $\mathbb{R}_{\geq0}$ are in one-to-one correspondence
with the continuous homomorphisms of $G$ to $\mathbb{R}_{>0}$. So
clearly $G$ admits a universal irreducible conic representation if
and only if it admits no non-trivial continuous homomorphisms to $\mathbb{R}_{>0}$.\\

\textbf{Example 7.1:} Compact Hausdorff groups can be shown to always
admit a ray, in any conic representation, which is not just invariant
but pointwise fixed. It is done in a strictly analogous manner to
the proof of their amenability, except one uses now the resultant
mapping instead of the barycenter (both are equivariant). So the only
irreducible representation of a compact group is the identity action
on $\mathbb{R}_{\geq0}$.\\

\textbf{Theorem 25:} A group $G$ that is amenable but not Tychonoff
does not admit a universal irreducible conic representation.\\

\textbf{Proof:} If there is a universal irreducible conic representation
then it is degenerate - it has an invariant section on which the restriction
of the action is affine. So from amenability the universal irreducible
conic representation is just $\mathbb{R}_{+}$with $G$ acting as
the identity. However, $G$ is not Tychonoff, hence there exists a
conic representation of $G$ without an invariant ray. By Zorn's lemma,
it contains an irreducible sub-representation without an invariant
ray, in particular with sections consisting of more than one point,
and as such it can not be a factor of the universal irreducible conic
representation. Contradiction. $\blacksquare$\\

The group $E$ of transformations of the plane that is generated by
rotations and translations is the semi-direct product of the rotations
sub-group and the normal translations sub-group. The rotations sub-group
$K$ is just the circle, so it is compact and hence does not admit
non-trivial continuous homomorphisms to the reals. The conjugacy classes
of the translations sub-group $T\cong\mathbb{R}^{2}$ are the circles
around the origin and hence it also does not admit non-trivial continuous
homomorphisms to the reals. $E=KT$ and hence it also does not admit
non-trivial continuous homomorphisms to the reals. It is solvable
and thus amenable, but is known not to be Tychonoff (see \cite{key-2})
\footnote{Another solvable group which is not Tychonoff is the one considered
in example 3.1. In that example the action has no invariant ray.}. From prop. 24, we also know it does not admit a universal irreducible
conic representation, and thus we conclude that not having non-trivial
continuous homomorphisms to the reals is not a sufficient condition
for having a universal irreducible conic representation. As we shall
see, this is the situation in the case of $SL_{2}\left(\mathbb{R}\right)$.\\

\section{Semi-Conic Representations}

We define a $semi-cone$ to be a sub-set of a real or complex linear
space closed just under multiplication in non-negative scalars \footnote{Some authors use the term $convex\,cone$ for $cone$, and $cone$
for $semi-cone$. }. For a semi-cone $W$ lying in a topological linear space with a
point-separating continuous dual, we define a section to be $L^{-1}\left(1\right)$
of a continuous homogeneous function $L:W\rightarrow\mathbb{\mathbb{R}}_{\geq0}$
that is positive on $W\setminus0$. For $W$ that admits a compact
section, a semi-conic representation of a topological group $G$ on
it is a continuous action $\tau:G\times W\rightarrow W$ such that
$\tau_{g}\left(\cdot\right)$ is homogeneous for any fixed $g\in G$.
A semi-conic representation of $G$ induces in a natural way an action
of $G$ on every section of the semi-cone together with a multiplier
of $G$ on the section. Conversely, given such a pair - a.k.a. $semi-conic\,pair$
- it induces a semi-conic representation on the semi-cone \footnote{If it is given without a cone, notice that any compact Hausdorff space
$Y$ is embeddable in a locally convex topological linear space $X$
(As the extreme points of $\Pr\left(Y\right)$), and thus $Y$ is
a section of the semi-cone in $X\varoplus\mathbb{R}$ generated by
$\left(Y,1\right)$. The contiuous dual of a locally convex topological
linear space separates points.}. Induced actions on different sections are naturally isomorphic as
topological dynamical systems. Homomorphisms of semi-conic representations
are also defined in a straightforward manner, and satisfy properties
analogous to the ones described in prop. 16 for conic representations.\\

A semi-conic representation of a group on a semi-cone is called $minimal$
if it has no sub-representations other than itself and $\left\{ 0\right\} $.
That is equivalent to the induced actions on the sections being mininal
in the category of topological dynamical systems. \\

A conic representation of $G$ on $V$ induces in a natural way a
semi-conic representation of $G$ on $\overline{ExtV}$, where $ExtV$
is the union of the extreme rays of $V$ and $\overline{ExtV}$ is
its closure (Not to be confused with $Ext$ when taken on a CCS).
As mentioned in the previous section and will be proven later, there
exists a reducible conic representation of $SL_{2}\left(\mathbb{R}\right)$
with a strongly proximal action on $\overline{Ext}$ of its sections.
The same cone admits an irreducible degenerate conic representation
with an identical action on $\overline{Ext}$ of its sections. This
means that given a conic representation $V$ of a group $G$ together
with a section $Q$, the induced action of $G$ on $\overline{ExtQ}$
does not determine whether $V$ is irreducible or not (unlike the
induced action on $Q$). However, we shall see later in this section
that the induced action of $G$ on $\overline{ExtV}$ does determine
this (prop. 26). \\

As in affine representations of $G$, $\overline{Ext}$ is generally
not a functor. If however we restrict to the category of irreducible
conic representations of $G$, and $\varphi:V_{1}\rightarrow V_{2}$
is a homomorphism of such, then by lemma 17 (to be proved), $\overline{ExtV_{1}}$
and $\overline{ExtV_{2}}$ are minimal and $\varphi\left(\overline{ExtV_{1}}\right)\supseteq\overline{ExtV_{2}}$,
hence $\varphi\left(\overline{ExtV_{1}}\right)=\overline{ExtV_{2}}$.
So $\overline{Ext}_{*}$ is a functor between that category and the
category of minimal semi-conic representations of $G$. In the opposite
direction, we do not need to restrict ourselves, and we have a functor
between the category of semi-conic representations of $G$ to the
category of its conic representations, which we will now introduce.\\

Let $W$ be a semi-conic representation of $G$ with action $\tau:G\times W\rightarrow W$.
Take $Y$ to be some section of $W$ defined by a conic function $L$,
And take $MSec\left(W\right)$ (abbr. of ``Measures on a Section'')
to be $M\left(Y\right)$ - the cone of non-negative regular borel
measures on $Y$ (with the weak-{*} topology as usual). Define a semi-conic
embedding $\varphi:W\rightarrow MSec\left(W\right)$ by $\varphi\left(z\right)=L\left(z\right)\delta_{\frac{z}{L\left(z\right)}}$,
and from now on $W$ is to be identified with its image under this
embedding. $G$ has an action on $\varphi\left(W\right)$ induced
by this identification which we will denote by $\tilde{\tau}$. We
want to extend the action $\tilde{\tau}$ of $G$ from $\varphi\left(W\right)$
to a conic action on all $MSec\left(W\right)$. We thus define the
action of $g\in G$ on $\mu\in MSec\left(W\right)=M\left(Y\right)$
with the aid of the resultant function $r:M_{C}\left(MSec\left(W\right)\right)\rightarrow MSec\left(W\right)$
to be $\tilde{\tau}_{g}\left(\mu\right)=r\left(\tilde{\tau}_{g*}\varphi_{*}\mu\right)$.
The meaning of the asteriks being, as usual, the ordinary push-forward
of measures. $\tilde{\tau}_{g*}$ being already well defined on measures
supported on $\varphi\left(W\right)$. The process is illustrated
in Figure 1. \\

\begin{center}
\begin{figure}[H]
\begin{centering}
\subfloat[This is $\mu$. The section is $\Pr\left(Y\right)$.]{\begin{centering}
\includegraphics[scale=0.8]{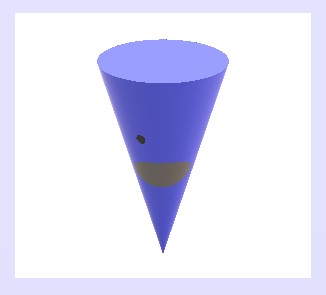}
\par\end{centering}
}~~\subfloat[This is the finite positive measure on $Y\protect\cong\overline{Ext\left(\Pr\left(Y\right)\right)}$
that $\mu$ represents. Notice that in the previous picture $\mu$
is not on $\Pr\left(Y\right)$ and thus in such a case it is not a
probability measure.]{\begin{centering}
\includegraphics[scale=0.8]{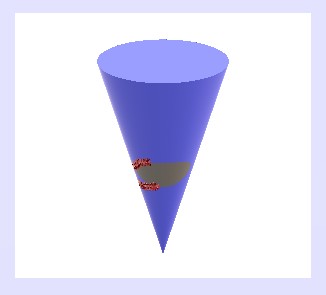}
\par\end{centering}

}
\par\end{centering}
\centering{}\subfloat[This is the push-forward of the measure in the previous picture through
the pre-defined action of $g$ on the semi-cone surrounding the cone
$M\left(Y\right)$ (this semi-cone is isomorphic to $W)$. ]{\begin{centering}
\includegraphics[scale=0.8]{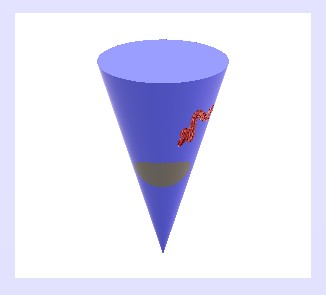}
\par\end{centering}
}~~\subfloat[This is the resultant in $M\left(Y\right)$ of the measure in the
previous picture.]{\centering{}\includegraphics[scale=0.8]{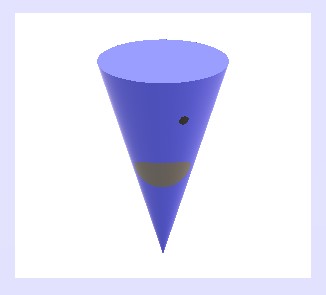}}\caption{A description of the definition of the action of $g\in G$ on a point
$\mu\in M\left(Y\right)$. The action of $g$ transfers the point
depicted in (a) to the point depicted in (d). }
\end{figure}
\par\end{center}

This definition of the action satisfies the algebraic requirements
for being a conic action since $r$ is equivariant and preserves non-negative
linear combinations. It remains to explain why it is continuous. For
any $0<a<b$ and a cone $V$ with a compact section $Q$, the resultant
mapping is continuous when restricted to the cone of the measures
that are supported on $[a,b]Q$. Any $g\in G$ has an open neighborhood
$U$ for which there exist $0<a<b$ such that $UY\subseteq[a,b]Y$.
Thus the action's restriction to $U\times MSec\left(W\right)$ is
continuous as a composition of continuous mappings, and thus it is
continuous. \\

It can be easily verified morphisms also transform as necessary, and
the construction of the covariant functor $MSec$ from semi-conic
representations to conic ones is done\footnote{It is easily verified to be independent of the choice of the section
up to a natural isomorphism. This justifies our preference of the
term $MSec\left(W\right)$ to $M\left(Y\right)$.}. In section 10 we will present an equivalent construction for $MSec$
that is somewhat longer but probably easier to digest.\\

The first part of the following proposition is virtually the purpose
for which $MSec$ was designed. Its second part implies that whether
a conic representation $V$ of $G$ is irreducible or not is determined
by the semi-conic representation $\overline{ExtV}$. We shall call
such a semi-conic representation an $irreducible\,semi-conic\,representation$
and its corresponding semi-conic pairs $irreducible\,semi-conic\,pairs$.\\

\textbf{Proposition 26:} Let $V$ be a conic representation of $G$
with action $\tau$, then $V$ is a factor of $MSec\left(V\right)$\footnote{$V$ is a cone, and a cone is in particular a semi-cone.}
through the resultant mapping $r$ (the same is true for $MSec\left(\overline{ExtV}\right)$
since it is embedded in $MSec\left(V\right)$). In addition, $V$
is irreducible if and only if $MSec\left(\overline{ExtV}\right)$
is irreducible.\\

\textbf{Proof:} Consider the section $Q$ of $V$ used in defining
$MSec\left(V\right)$ (it is $M\left(Q\right)$). As was already mentioned
$r:M_{C}\left(V\right)\rightarrow V$ preserves non-negative linear
combinations, and its restriction $r:M\left(Q\right)\rightarrow V$
is continuous since $M\left(Q\right)$ is a sub-cone of $M_{C}\left(V\right)$
that admits a compact section. For any $g\in G$ and $x\in Q$, $r\left(\tilde{\tau_{g}}\left(\delta_{x}\right)\right)=r\left(r\left(\tilde{\tau}_{g*}\varphi_{*}\delta_{x}\right)\right)=\tau_{g}\left(r\left(\delta_{x}\right)\right)$.
The mappings are linear and thus the equality also holds for non-negative
linear combinations of dirac measures, and are continuous and thus
- by Krein-Milman - we obtain $r\left(\tilde{\tau_{g}}\left(\mu\right)\right)=\tau_{g}\left(r\left(\mu\right)\right)$
for any $\mu\in M\left(Y\right)$. \\

For the second part, the 'if'' part is obvious since the pre-image
of a conic sub-representation through a conic homomorphism (the resultant)
is a sub-representation. For the ``only if'' part, take a section
$Q$ of $V$, and construct $MSec\left(\overline{ExtV}\right)$ using
the section $\overline{ExtQ}$ of $\overline{ExtV}$. $r|_{MSec\left(\overline{ExtV}\right)}^{-1}\left(Q\right)$
is the section $\Pr\left(\overline{ExtQ}\right)$, and $r|_{MSec\left(\overline{ExtV}\right)}:\Pr$$\left(\overline{ExtQ}\right)\rightarrow Q$
is in fact the barycentric mapping $\beta$, and it is a homomorphism
of the dynamical systems $\Pr\left(\overline{ExtQ}\right)$ and $Q$
under the induced actions on them which we denote by $\tilde{\rho}$
and $\rho$ respectively . If $\mu\in\Pr\left(\overline{ExtQ}\right)$
and $y\in\overline{ExtQ}$ then there exists a net $g_{\alpha}$ such
that $\rho_{g_{\alpha}}\left(\beta\left(\mu\right)\right)$ converges
to $y$ (by lemma 17) and $\tilde{\rho}_{g_{\alpha}}\left(\mu\right)$
converging to some $\nu\in\Pr\left(\overline{ExtQ}\right)$ (by the
compactness of $\Pr\left(\overline{ExtQ}\right)$). Since $\beta$
is a continuous homomorphism we have $\beta\left(\nu\right)=y$, hence
$\nu=\delta_{y}$. $\blacksquare$\\

\textbf{Proposition 27:} Let $W$ be a semi-cone with section $Y$,
and $MSec\left(W\right)=M\left(Y\right)$. Let $Q$ be a section of
$MSec\left(W\right)$ defined by the continuous conic funcion $L:MSec\left(W\right)\rightarrow\mathbb{R}_{\geq0}$
which is strictly positive on $MSec\left(W\right)\setminus0$ (i.e.
$Q=L^{-1}\left(1\right)$). Then there exists a positive function
$f\in C\left(Y\right)$ such that $L\left(\mu\right)=\int_{Y}f\left(y\right)d\mu\left(y\right)$.\\

\textbf{Proof:} If there exists such a function and $a\delta_{y}$
is an extreme point of $Q$ for some $a>0$, then $f\left(y\right)=\frac{1}{a}$.
We now have $f$ defined, it is continuous and we can use Krein-Milman
theorem to prove it is in fact equal to $L$. $\blacksquare$\\

So we reduced in $MSec\left(W\right)$ the section definition from
one using conic functions to one using linear functionals.\\

Assume now we have a section $Q_{1}$ in a conic representation $V$
of $G$ defined by the positive conic function $L_{1}$ on $V$. Given
another positive conic function $L_{2}$ on $V$, it defines another
section $Q_{2}$, and we may wonder how the multiplier $\sigma_{1}$
of $Q_{1}$ relates to the multiplier $\sigma_{2}$ of $Q_{2}$. \\

\textbf{Proposition 28:} $\sigma_{2}\left(g,\frac{y}{L_{2}\left(y\right)}\right)=\frac{L_{2}\left(\rho_{g}\left(y\right)\right)}{L_{2}\left(y\right)}\sigma_{1}\left(g,y\right)$
for all $g\in G$, $y\in Q_{1}$, where $\rho_{g}$ is the induced
action of $g$ on the section $Q_{1}$.

\textbf{Proof:} The statement is equivalent to $L_{2}\left(y\right)L_{2}\left(g\frac{y}{L_{2}\left(y\right)}\right)=L_{2}\left(\frac{gy}{L_{1\left(gy\right)}}\right)L_{1}\left(gy\right)$.
$\blacksquare$\\

An equivalent strictly analogous proposition can be stated for a semi-conic
representation $W$, its section $Y$ , and a continuous positive
and positive-homogeneous function $L$ on $W$ defining another section
$Y'$. The restriction of $L$ to $Y$ determines a continuous function
on $Y$, and vice versa, any continuous function on $Y$ can be extended
to a unique continuous positive and positive-homogeneous function
on $W$. So giving either of the two is essentially the same. Given
a continuous function $f$ on $Y$ we thus get the following corollary
relating the multipliers $\sigma_{1}$ of $Y$ and $\sigma_{2}$ of
the section defined by $f$.\\

\textbf{Corollary 29:} $\sigma_{2}\left(g,\frac{y}{f\left(y\right)}\right)=\frac{f\left(\rho_{g}\left(y\right)\right)}{f\left(y\right)}\sigma_{1}\left(g,y\right)$
for all $g\in G$, $y\in Y$, where $\rho_{g}$ is the induced action
of $g$ on the section $Y$.\\

We develop our jargon a bit futher. If $\left(Y,\rho\right)$ is a
topological dynamical system of $G$, and $h\in C\left(Y\right)$
is positive, then $\frac{h\left(\rho_{g}\left(y\right)\right)}{h\left(y\right)}$
is a multiplier of $G$ on $Y$. We call the multipliers of this form
$trivial\,multipliers$. The trivial multipliers form a multiplicative
group that we denote by $T\left(Y\right)$. We denote by $M\left(Y\right)$
the group of all multipliers of $G$ on $Y$ and by $H\left(Y\right)$
the quotient group $M\left(Y\right)/T\left(Y\right)$. Stating our
previous result in this new terminology we obtain\\

\textbf{Theorem 30:} Given a semi-conic pair $\left(\rho,\sigma\right)$
of $G$ on $Y$, the semi-conic pairs of all of the sections of the
semi-conic representation it induces are exactly all pairs of the
form $\left(\rho,\sigma'\right)$ on $Y$ where $\sigma'$ belongs
to the coset $\sigma T\left(Y\right)$ (under the canonical identification
between sections).\\

\textbf{Corollary 31:} Given a conic representation $V$ of $G$ and
one of its sections $Q$, $V$ is degenerate if and only if the multiplier
on $\overline{ExtQ}$ is trivial.\\

\textbf{Corollary 32:} Given two semi-conic pairs $\left(\rho,\sigma_{1}\right)$,
$\left(\rho,\sigma_{2}\right)$ of $G$ on $Y$. The two semi-conic
representations induced by the pairs are isomorphic if and only if
$\sigma_{2}\in\sigma_{1}T\left(Y\right)$.\\

\textbf{Lemma 33:} If $Y$ is a compact Hausdorff space, $G$ acts
continuously on $Y$ by $\rho:G\times Y\rightarrow Y$ and the action
is minimal. Then any \textbf{bounded} multiplier $\sigma:G\times Y\rightarrow\mathbb{R}_{>0}$
is trivial.\\

\textbf{Proof:} In our terminology, $\rho$ and $\sigma$ form a semi-conic
pair. $Y$ is embeddable in a locally convex linear space $X$, and
we now consider the semi-conic representation $\tau$ the pair induces
on the semi-cone in $X\varoplus\mathbb{R}$ generated by $\left(Y,1\right)$.
We denote the projection onto the second summand by $L:X\varoplus\mathbb{R\rightarrow\mathbb{R}}$.\\

Since $\sigma$ is bounded and $Y$ is compact the invariant set $\overline{\tau_{G}\left(Y,1\right)}$
is compact, hence it contains a minimal set $Z$ (in the category
of ordinary topological dynamics). Defining $\varphi:Z\rightarrow\left(Y,1\right)$
by $\varphi\left(z\right)=\frac{z}{L\left(z\right)}$. $\varphi$
is onto since when the action on $\left(Y,1\right)$ is taken to be
$\rho$ (under its natural identication with $Y$) which is minimal,
it commutes with the actions.\\

We claim $\varphi$ is also one-to-one, and therefore $\varphi$ is
invertible (the inverse in continuous because $\varphi$ is a closed
map). If it were not, there would exist $z,z'\in Z$ such that $z'=rz$
for some $r>1$. This implies the value of $L$ on $Z$ is unbounded,
in contradiction to $Z$ being compact. For if $L\left(\tau_{g_{1}}z\right)=M$
then $L\left(\tau_{g_{1}}\left(z'\right)\right)=L\left(\tau_{g_{1}}\left(rz\right)\right)=rM$,
and hence for any $\epsilon>0$ there is a neighborhood $U_{\epsilon}$
of $z'$ such that $L\left(\tau_{g_{1}}\left(U_{\epsilon}\right)\right)\subseteq(rM-\epsilon,\infty)$.
Since $Z$ is minimal there exists $g_{2}\in G$ for which $\tau_{g_{2}}\left(z\right)\in U_{\epsilon}$.
Thus $L\left(\tau_{g_{1}g_{2}}\left(z\right)\right)>rM-\epsilon$
and the value of $L$ on $Z$ is unbounded. \\

Defining $f:Y\rightarrow\mathbb{R}_{>0}$ by $f\left(y\right)=L\left(\varphi^{-1}\left(y,1\right)\right)$,
we obtain $\sigma\left(g,y\right)=\frac{f\left(\rho_{g}\left(y\right)\right)}{f\left(y\right)}$.
$\blacksquare$\\

\textbf{Theorem 34:} An irreducible factor $V$ of a degenerate conic
representation of $G$ is itself degenerate.\\

\textbf{Proof:} Taking $Q$ a section of $V$ with induced conic pair
$\left(\rho,\sigma\right)$, $\sigma$ is bounded by lemma 22. Thus
so is its restriction to the minimal set $\overline{ExtQ}$ with respect
to the action $\rho$ (lemma 17). By lemma 33, that restriction is
a trivial multiplier, and by cor. 31 we are done. $\blacksquare$\\

On the other hand we have the following example.\\

\textbf{Example 8.1:} We now fulfill an obligation from the past (see
the remark after lemma 22), and present an example showing that a
non-irreducible factor of a degenerate conic representation need not
be degenerate (although by lemma 22 the multipliers of its sections
are bounded). Using the results obtained in this section ($MSec$
and cor. 31) it is sufficient to construct a compact Hausdorff $G$-space
$Y$, a factor $S$ of $Y$, and a non-trivial multiplier on $S$
such that its pullback to a multiplier on $Y$ is trivial. For $G$
we take $\mathbb{Z}$, and for $Y$ we take $\left\{ 0,1\right\} \cup\left\{ a_{n}\right\} _{n\in\mathbb{Z}}$
where $\left\{ a_{n}\right\} _{n\in\mathbb{Z}}\subseteq\left(0,1\right)$
is any sequence satisfying $a_{n+1}<a_{n}$ for all $n\in\mathbb{Z}$
and $lim_{n\rightarrow\infty}a_{n}=0$, $lim_{n\rightarrow-\infty}a_{n}=1$.
We define $T:Y\rightarrow Y$ by \\

\[
T\left(a_{n}\right)=a_{n+1},\,T\left(0\right)=0,\,T\left(1\right)=1.
\]
\\
The action of $\mathbb{Z}$ is defined by $T$ ($1$ acts on $Y$
as $T$ ). To define a multiplier $\sigma:\mathbb{Z}\times Y\rightarrow\mathbb{R}_{>0}$
on $Y$ we take any positive sequence $\left\{ r_{n}\right\} _{n\in\mathbb{\mathbb{\mathbb{Z}}}}$
such that $\prod_{n=1}^{\infty}r_{n}=\frac{1}{2}$ and $r_{n}=1$
for any non-positive $n$. We define $\sigma\left(k,0\right)=\sigma\left(k,1\right)=\sigma\left(0,y\right)=1$
for $k\in\mathbb{Z}$ and $y\in Y$ , and $\sigma\left(k,a_{n}\right)=r_{n+1}\cdot...\cdot r_{n+k}$
for $k\geq1$, $\sigma\left(k,a_{n}\right)=r_{n}^{-1}\cdot...\cdot r_{n+k+1}^{-1}$
for $k<0$. To obtain $S$ we identify $0$ and $1$ and denote the
quotient mapping by $\pi:Y\rightarrow S$. Notice that $S$ is Hausdorff
and that the quotient respects the action of $\mathbb{Z}$ on $Y$
and the multiplier $\sigma$. We shall denote the resulting multiplier
on $S$ by $\eta$ ($\sigma\left(k,y\right)=\eta\left(k,\pi\left(y\right)\right)$).\\

We claim $\eta$ is non-trivial. Assuming the contrary, there exists
a continuous function $f:S\rightarrow\mathbb{R}_{>0}$ such that $\eta\left(k,y\right)=\frac{f\left(T^{k}y\right)}{f\left(y\right)}$
and we assume without loss of generality that $f\left(\pi\left(a_{0}\right)\right)=1$.
So $f\left(\pi\left(a_{n}\right)\right)=f\left(\pi\left(T^{n}a_{0}\right)\right)=\eta\left(n,\pi\left(a_{0}\right)\right)$,
and this is equal for $n\geq1$ to $\prod_{k=1}^{n}r_{k}$ and for
$n<0$ to $1$. Thus $f\left(\pi\left(0\right)\right)=\prod_{n=1}^{\infty}r_{n}=\frac{1}{2}$
and $f\left(\pi\left(1\right)\right)=1$, but $\pi\left(0\right)=\pi\left(1\right)$
and that is a contradiction, hence $\eta$ is non-trivial. From this
reasoning it is also clear that $\sigma$ is trivial (one just defines
$f:Y\rightarrow\mathbb{R}_{>0}$ by these requirements), and we are
done.\\

\textbf{Corollary 35:} $G$ admits a universal irreducible conic representation
if and only if all irreducible conic representations of $G$ are degenerate.\\

\textbf{Proof:} The degenerate conic representation induced from the
universal affine representation of $G$ admits a homomorphism to any
other irreducible degenerate conic representation by extending the
corresponding homomorphism of affine systems. The ``only if'' part
is a direct consequence of prop. 20 and theorem. 34. Another way to
obtain this result is by recalling that if $Q$ is a section of the
the universal irreducible conic representation, the induced action
on $\overline{ExtQ}$ is proximal (strongly proximal) and then use
lemma 36. $\blacksquare$\\

\textbf{Lemma 36:} Let $Y$ and $Z$ be compact Hausdorff spaces with
$G$ acting on them continuously. Let $\varphi:Y\rightarrow Z$ be
a continuous equivariant mapping which is onto, $\sigma_{Z}$ a multiplier
on $Z$, and $\sigma_{Y}$ is its pull back to $Y$ for any $g\in G$.
If $\sigma_{Y}$ is a trivial multiplier on $Y$, and the action on
$Y$ is \textbf{proximal} then $\sigma_{Z}$ is also trivial.\\

\textbf{Proof:} Denoting the action on $Y$ by $\rho:G\times Y\rightarrow Y$,
we have $\sigma_{Y}\left(g,\varphi\left(x\right)\right)=\frac{h\left(\rho_{g}\left(x\right)\right)}{h\left(x\right)}$
for some positive $h\in C\left(Y\right)$. We now show $h$ respects
the fibers of $\varphi$ and this ends the proof since $\varphi$
is a quotient mapping.\\

Let $x,y\in Y$ such that $\varphi\left(x\right)=\varphi\left(y\right)$.
That is to all $g\in G$\\

$\frac{h\left(\rho_{g}\left(x\right)\right)}{h\left(x\right)}=\frac{h\left(\rho_{g}\left(y\right)\right)}{h\left(y\right)}$\\

$\frac{h\left(\rho_{g}\left(y\right)\right)}{h\left(\rho_{g}\left(x\right)\right)}=\frac{h\left(y\right)}{h\left(x\right)}$\\

But $Y$ is proximal and hence $\frac{h\left(y\right)}{h\left(x\right)}=1$,
that is $h\left(x\right)=h\left(y\right)$. $\blacksquare$

\section{Prime Conic Systems}

A conic representation of a topological group $G$ is called $prime$
if all its homomorphisms to other representations of $G$ - that are
not one-dimensional - are one-to-one. A Semi-conic representation
$W$ of $G$ is called $Conically\,Prime$ if $MSec\left(W\right)$
is prime.\\

Let $Y$ be a topological dynamical system of $G$. $Y$ is called
$affinely\,prime$ if $\Pr\left(Y\right)$ is a prime affine dynamical
system (meaning all its homomorphisms to non-trivial representations
are one-to-one). Given $f\in C\left(Y\right)$\footnote{The space of continuous real valued functions on $Y$.}
which is not constant we define $V_{f}=\overline{Span\left\{ f\circ\rho_{g}:g\in G\right\} }$
(the closure taken in the uniform norm). The system is said to have
the $Linear\,Stone-Weierstrass$ (LSW) property if for any such $f$,
the direct sum of $V_{f}$ and the space of constant functions is
all $C\left(Y\right)$. It is known that being LSW is equivalent to
being affinely prime (to be found in \cite{key-5}). We will not use
this fact, but in a strictly analogous proof we will show the following
proposition.\\

\textbf{Propositon 37:} Let $W$ be a semi-conic representation of
$G$ and let $Y$ be its section. If the induced action of $G$ on
$Y$ is LSW then $W$ is conically prime. \\

\textbf{Proof:} Let $\pi:MSec\left(W\right)\rightarrow V$ be a homomorphism
of $MSec\left(W\right)$ to some conic representation $V$. Given
$f\in C\left(Y\right)$, we denote by $\hat{f}$ the conic extenstion
of $f$ to all $MSec\left(W\right)$, i.e. $\hat{f}\left(\mu\right)=_{Y}\int f\left(y\right)d\mu\left(y\right)$.
Now if $F$ is a continuous conic function on $V$, then its pullback
$F\circ\pi$ is a continuous conic function on $MSec\left(W\right)$.
The set \\
\[
\left\{ f\in C\left(Y\right):\hat{f}\,is\,such\,a\,pullback\,through\,\pi\right\} 
\]

is a closed sub-space of $C\left(Y\right)$ invariant under the action
of $G$, contains a non-constant function and all constant functions,
therefore - by the LSW property - it is equal to all $C\left(Y\right)$.
If we have $\mu,\kappa\in MSec\left(W\right)$ such that $\pi\left(\mu\right)=\pi\left(\kappa\right)$,
then all pull-backs as above are equal on $\mu$ and $\kappa$, so
$\hat{f}\left(\mu\right)=\hat{f}\left(\kappa\right)$ for all $f\in C\left(Y\right)$,
which is by definition $\int f\left(y\right)d\mu\left(y\right)=\int f\left(y\right)d\kappa\left(y\right)$,
which means $\mu=\kappa$. $\blacksquare$\\

\textbf{Example 9.1:} The universal strongly proximal topological
dynamical system of the group $SL_{2}\left(\mathbb{R}\right)$ in
known to be LSW (the proof can be found in \cite{key-5}), hence its
degenerate conic representation generated by its unique irreducible
affine representation is prime. However, if $SL_{2}\left(\mathbb{R}\right)$
admits a universal irreducible conic representation it is the latter,
and so if this is the case then it is its unique irreducible conic
representation. However, we will see it is not unique.\\

\section{An Alternative Approach to Construct the Group Action on MSec$\left(W\right)$}

Given a semi conic representation of $G$ on a semi-cone $W$. Our
approach will now be to take a section $Y$ and extend the semi-conic
pair $\left(\rho,\sigma\right)$ on $Y$ to a conic pair $\left(\tilde{\rho},\tilde{\sigma}\right)$
on $\Pr\left(Y\right)$. Thus getting a conic action of $G$ on on
$MSec\left(W\right)$. It is equivalent to our original definition
since the construction presented will be easily seen to be the unique
extension of the semi-conic pair $Y$ to a conic pair on $\Pr\left(Y\right)$.\\

Since fixing any $g\in G$, $\tilde{\sigma}\left(g,\cdot\right)$
should be a continuous affine function on $\Pr\left(Y\right)$, then
if it exists it necessarily satisfies for every $\text{\ensuremath{\nu\in\Pr\left(Y\right)}}$:
$\tilde{\sigma}\left(g,\nu\right)=\int_{Y}\sigma\left(g,y\right)d\nu\left(y\right)$
($\nu$ is the barycenter of itself). By identification of the measures
with linear functionals, and recalling the definition of the weak-{*}
topology, one is easily convinced that the above formula indeed defines
a continuous $\tilde{\sigma}$ on $\Pr\left(Y\right)$. It is also
obviously affine. The only thing left for checking is that it is a
multiplier, but first we should extend the definition of $\rho_{g}$
to all $\Pr\left(Y\right)$: \\
\[
\tilde{\rho_{g}}\left(\nu\right)=\frac{\int_{Y}\sigma\left(g,y\right)\delta_{\rho_{g}\left(y\right)}d\nu\left(y\right)}{\tilde{\sigma}\left(g,\nu\right)}
\]
\\
As in the definition of the barycenter and resultant, the integration
here is of Pettis kind, and it works for similar reasons. This definiton
of $\rho_{g}$ was conceived just in order for it to satisfy $\left(1\right)$
(in the proof of prop. 10) so it is no surprise that it does. As promised,
we are ready to check now that $\sigma$ is a multiplier on all $\Pr\left(Y\right)$:\\

On the one hand, by definition\\
\[
\tilde{\sigma}\left(g\gamma,\nu\right)=\int_{Y}\sigma\left(g\gamma,y\right)d\nu\left(y\right)
\]
But on the other,

\[
\tilde{\sigma}\left(\gamma,\nu\right)\cdot\tilde{\sigma}\left(g,\tilde{\rho_{\gamma}}\left(\nu\right)\right)=\tilde{\sigma}\left(\gamma,\nu\right)\int_{Y}\sigma\left(g,u\right)d\left(\tilde{\rho_{\gamma}}\left(\nu\right)\right)\left(u\right)=\tilde{\sigma}\left(\gamma,\nu\right)\cdot\frac{1}{\tilde{\sigma}\left(\gamma,\nu\right)}\int_{Y}\sigma\left(\gamma,y\right)\sigma\left(g,\rho_{\gamma}\left(y\right)\right)d\nu\left(y\right)
\]
\\
The last equality follows by considering $\nu=\lambda_{1}\delta_{y_{1}}+...+\lambda_{n}\delta_{y_{n}}$,
a convex combination of Dirac measures (which is dense in $\Pr\left(Y\right)$
by the Krein-Milman theorem): $\tilde{\rho_{\gamma}}\text{\ensuremath{\left(\nu\right)}}=\frac{\int_{Y}\sigma\left(\gamma,y\right)\delta_{\rho_{\gamma}\left(y\right)}d\nu\left(y\right)}{\tilde{\sigma}\left(\gamma,\nu\right)}=\frac{\sum_{i=1}^{n}\lambda_{i}\sigma\left(\gamma,y_{i}\right)\delta_{\rho_{\gamma}\left(y_{i}\right)}}{\tilde{\sigma}\left(\gamma,\nu\right)}$
and\\

$\int_{Y}\sigma\left(g,u\right)d\left(\tilde{\rho_{\gamma}}\left(\nu\right)\right)\left(u\right)=\frac{1}{\tilde{\sigma}\left(\gamma,\nu\right)}\sum_{i=1}^{n}\int_{Y}\lambda_{i}\sigma\left(\gamma,y_{i}\right)\sigma\left(g,u\right)d\left(\delta_{\rho_{\gamma}\left(y_{i}\right)}\right)\left(u\right)=\frac{1}{\tilde{\sigma}\left(\gamma,\nu\right)}\sum_{i=1}^{n}\int_{Y}\lambda_{i}\sigma\left(\gamma,y_{i}\right)\sigma\left(g,\rho_{\gamma}\left(y\right)\right)d\left(\delta_{y_{i}}\right)\left(y\right)=$\\

$=\frac{1}{\tilde{\sigma}\left(\gamma,\nu\right)}\sum_{i=1}^{n}\int_{Y}\sigma\left(\gamma,y\right)\sigma\left(g,\rho_{\gamma}\left(y\right)\right)d\left(\lambda_{i}\delta_{y_{i}}\right)\left(y\right)=\frac{1}{\tilde{\sigma}\left(\gamma,\nu\right)}\sum_{i=1}^{n}\int_{Y}\sigma\left(\gamma,y\right)\sigma\left(g,\rho_{\gamma}\left(y\right)\right)d\left(\lambda_{i}\delta_{y_{i}}\right)\left(y\right)=\frac{1}{\tilde{\sigma}\left(\gamma,\nu\right)}\int_{Y}\sigma\left(\gamma,y\right)\sigma\left(g,\rho_{\gamma}\left(y\right)\right)d\nu\left(y\right)$.\\

And so, using the multiplier property of $\sigma$, we obtain $\tilde{\sigma}\left(\gamma,\nu\right)\cdot\tilde{\sigma}\left(g,\rho_{\gamma}\left(\nu\right)\right)=\tilde{\sigma}\left(g\gamma,\nu\right)$
and we have shown that $\left(\tilde{\rho},\tilde{\sigma}\right)$
is a conic pair on $\Pr\left(Y\right)$.\\

As already mentioned, the purpose of the whole construction of $MSec$
was the first part of prop. 26. In the terminology of this section
it is equivalent to the statement that given a $CCS$ $Q$ with a
conic pair $\left(\rho,\sigma\right)$, the barycenter mapping $\beta:\Pr\left(Q\right)\rightarrow Q$
commutes with $\tilde{\rho_{g}}$ and $\rho_{g}$ for any $g\in G$.
Because of the importance of this result let us give another independent
proof that it holds this time using the construction of this section.\\

First note that $\left(1\right)$ (in the proof of prop. 10) is equivalent
to requiring that any $x_{1},...,x_{n}\in Q$ and $\lambda_{1},...,\lambda_{n}\geq0$
such that $\lambda_{1}+...+\lambda_{n}=1$ satisfy: \\

$\rho_{g}\left(\lambda_{1}x_{1}+...+\lambda_{n}x_{n}\right)=\frac{\lambda_{1}\sigma\left(g,x_{1}\right)}{\sigma\left(g,\lambda x_{1}+..+\lambda_{n}x_{n}\right)}\rho_{g}\left(x_{1}\right)+...+\frac{\lambda_{n}\sigma\left(g,x_{n}\right)}{\sigma\left(g,\lambda x_{1}+..+\lambda_{n}x_{n}\right)}\rho_{g}\left(x_{n}\right)$
because $\sigma\left(g,\cdot\right)$ is affine for any fixed $g\in G$.\\

This in turn implies, by the Krein-Milman theorem, that $\left(1\right)$
is equivalent to requiring that for any $\nu\in\Pr\left(Q\right)$:\\

$\rho_{g}\left(\beta\left(\nu\right)\right)=\int\frac{\sigma\left(g,x\right)}{\sigma\left(g,\beta\left(\nu\right)\right)}\rho_{g}\left(x\right)d\nu\left(x\right)$.\\

Hence, our conic pair satisfies $\rho_{g}\left(\beta\left(\nu\right)\right)=\int\frac{\sigma\left(g,x\right)}{\sigma\left(g,\beta\left(\nu\right)\right)}\rho_{g}\left(x\right)d\nu\left(x\right)=\int\frac{\sigma\left(g,x\right)}{\tilde{\sigma}\left(g,\nu\right)}\rho_{g}\left(x\right)d\nu\left(x\right)$.
On the other hand, $\beta\left(\tilde{\rho_{g}}\left(\nu\right)\right)=\beta\left(\frac{\int_{Y}\sigma\left(g,x\right)\delta_{\rho_{g}\left(x\right)}d\nu\left(x\right)}{\tilde{\sigma}\left(g,\nu\right)}\right)=\frac{\int_{Y}\sigma\left(g,x\right)\rho_{g}\left(x\right)d\nu\left(x\right)}{\tilde{\sigma}\left(g,\nu\right)}$,
and the independent proof is done.\\

We now calculate the Radon-Nikodym derivative of $\tilde{\rho_{g}}\left(\nu\right)$
with respect to $\left(\rho_{g}\right)_{*}\nu$.\\

\textbf{Proposition 38:} $d\left(\tilde{\rho_{g}}\left(\nu\right)\right)\left(y\right)=\frac{d\left(\rho_{g*}\nu\right)\left(y\right)}{\sigma\left(g^{-1},y\right)\sigma\left(g,\nu\right)}$
($\rho_{g*}\nu$ is the ordinary push forward of the measure $\nu$
by $\rho_{g}:Y\rightarrow Y$).\\

\textbf{Proof:} $\tilde{\rho_{g}}\left(\nu\right)=\frac{\int_{Y}\sigma\left(g,y\right)\rho_{g}\left(y\right)d\nu\left(y\right)}{\tilde{\sigma}\left(g,\nu\right)}=\int_{Y}\frac{\sigma\left(g,y\right)}{\tilde{\sigma}\left(g,v\right)}\rho_{g}\left(y\right)d\nu\left(y\right)=\int_{Y}\frac{\sigma\left(g,y\right)}{\tilde{\sigma}\left(g,v\right)}\delta_{gy}\left(y\right)d\nu\left(y\right)=\int_{Y}\frac{\sigma\left(g,g^{-1}y\right)}{\tilde{\sigma}\left(g,v\right)}\delta_{y}\left(y\right)d\left(\rho_{g*}\nu\right)\left(y\right)$
\\
Hence, for any $f\in C\left(Y\right)$, the integral of $f$ with
respect to $\tilde{\rho_{g}}\left(\nu\right)$ is , by the definintion
of the Pettis integral:\\

$\int_{Y}\frac{\sigma\left(g,g^{-1}y\right)}{\tilde{\sigma}\left(g,v\right)}f\left(y\right)d\left(\rho_{g*}\nu\right)\left(y\right)$.\\

Hence\\

$d\left(\tilde{\rho_{g}}\left(\nu\right)\right)\left(y\right)=\frac{\sigma\left(g,g^{-1}y\right)}{\tilde{\sigma}\left(g,\nu\right)}d\left(\rho_{g*}\nu\right)\left(y\right)=\frac{d\left(\rho_{g*}\nu\right)\left(y\right)}{\sigma\left(g^{-1},y\right)\tilde{\sigma}\left(g,\nu\right)}$~$\blacksquare$\\

\section{The Case Where $G$ is a Semi-Simple Lie Group}

We now consider the linear space of continuous functions on the associated
symmetric space $D$ of a semi-simple Lie group $G$, i.e. $D=K$\textbackslash$G${}
where $K$ is a maximal compact sub-group of $G$. This space can
be identified with the linear space of continuous functions $f$ on
$G$ satisfying $f\left(kg\right)=f\left(g\right)$ for any $g\in G$,
$k\in K$. Imposing the topology of pointwise convergence, $G$ acts
linearly and continuously on this locally convex linear space by translation
from the right. We denote by $V_{D}$ the cone of continuous positive
functions on $D$, but notice that it does not admit a compact section.\\

In \cite{key-2} the irreducible conic representations embedded in
$V_{D}$ are completely characterized. A $K-multiplier$ of an action
of $G$ on some compact Hausdorff space $Y$ is a multiplier $\sigma:G\times Y\rightarrow\mathbb{R}_{>0}$
such that $\sigma\left(k,y\right)=1$ for all $k\in K$, $y\in Y$.
Given a $K$-multiplier $\sigma$ on $Y$ we denote by $V\left(\sigma\right)$
the closed cone in $V_{D}$ generated by the set of functions $\left\{ \sigma\left(\cdot,y\right)|y\in Y\right\} $.
$V\left(\sigma\right)$ is closed under the action of $G$, and as
we shall see it admits a compact section, hence it is a conic representation
of $G$. In fact, it is shown in \cite{key-2} that the irreducible
conic representations in $V_{D}$ are exactly $V\left(\sigma\right)$
for $\sigma$'s which belong to a certain class of multipliers on
the universal strongly proximal space of the group called in \cite{key-2}
$irreducible$ $K-multipliers$ (we will not give their definition
here). \\

\textbf{Lemma 39:} Let $Y$ be a compact set in a topological linear
space with a point-separating continuous dual, then $\overline{Conv\left(Y\right)}$
is compact and $Ext\left(\overline{Conv\left(Y\right)}\right)\subseteq Y$.\\

\textbf{Proof:} The barycentric mapping $\beta:\Pr\left(Y\right)\rightarrow\overline{Conv\left(Y\right)}$
has an image which is a CCS containing $Y$ and thus it is onto and
$\overline{Conv\left(Y\right)}$ is compact. However, extreme points
are only barycenters of Dirac measures and hence all extreme points
of $\overline{Conv\left(K\right)}$ belong to $Y$. $\blacksquare$\\

\textbf{Corollary 40:} If a topological group $G$ acts linearly and
continuously on a topological linear space with a point-separating
continuous dual, and $Y$ is an invariant compact set on which the
action is transitive, then $Y=Ext\left(\overline{Conv\left(Y\right)}\right)$.\\

\textbf{Proof:} We need to prove $Y\subseteq Ext\left(\overline{Conv\left(Y\right)}\right)$.
$Ext\left(\overline{Conv\left(Y\right)}\right)$ is non-empty - say
by Krein-Milman - and hence the lemma implies $Y$ contains an extreme
points of $\overline{Conv\left(Y\right)}$. Since $G$ acts transitively
and preserves the affine structure, all points of $Y$ are extreme
points of $\overline{Conv\left(Y\right)}$. $\blacksquare$\\

Let $\sigma$ be a $K$-multiplier of $G$ on a compact Hausdorff
\textbf{homogeneous} space $Y$. $y\mapsto\sigma\left(\cdot,y\right)$
is a closed continuous mapping $\varphi:Y\rightarrow V_{D}$. The
evaluation of functions at the identity $e\in G$ is a continuous
linear functional, and thus $V\left(\sigma\right)$ has a compact
section $\overline{Conv\left(\varphi\left(Y\right)\right)}$ comprised
of its functions that satisfy $f\left(e\right)=1$. By cor. 40, $Ext\left(\overline{Conv\left(\varphi\left(Y\right)\right)}\right)=\varphi\left(Y\right)$.
If $Y$ is a section of a semi-conic representation $W$ with multiplier
$\sigma$, then $\varphi$ can be extended to $\varphi:W\rightarrow ExtV\left(\sigma\right)$.
Now $\left(\gamma\varphi\left(y\right)\right)\left(g\right)=\varphi\left(y\right)\left(g\gamma\right)=\sigma\left(g\gamma,y\right)=\sigma\left(\gamma,y\right)\sigma\left(g,\gamma y\right)=\sigma\left(\gamma,y\right)\cdot\varphi\left(\gamma y\right)\left(g\right)$
and hence by definiton, the multiplier $\eta:G\times\varphi\left(Y\right)\rightarrow\mathbb{R}_{>0}$
of the representation on $\varphi\left(Y\right)$ satisfies $\eta\left(\gamma,\varphi\left(y\right)\right)=\sigma\left(\gamma,y\right)$,
and $\varphi$ is homomorphism in the category of semi-conic representations
of $G$.\\

In a similar manner to the above, if $\left(\rho,\sigma\right)$ is
a conic (and not semi-conic) pair of $G$ on a CCS $Q$, and $\sigma$
is a $K$-multiplier, then $V\left(\sigma\right)$ is a factor of
the conic representation the pair induces (the mapping is affine since
$\sigma$ is affine in its second variable). For this we need not
even require the action of $G$ to be transitive on $\overline{ExtQ}$.
Notice that in this case, taking accordingly $\varphi:Q\rightarrow V_{D}$
with $x\mapsto\sigma\left(\cdot,x\right)$ for $x\in Q$, we have
$\varphi\left(Q\right)=\overline{Conv\left(\varphi\left(Q\right)\right)}$
since an affine image of a $CCS$ is a $CCS$.\\

By the above mentioned characterization of the irreducible conic reperesentations
in $V_{D}$ found in \cite{key-2}, we obtain that if $V$ is such
a representation, then the sections of $\overline{ExtV}$ are strongly
proximal. There is no example which is known to us of an irreducible
conic representation of a semi-simple group which does not have strongly
proximal sections. For reasons that will become apparent later, maybe
$V_{D}$ should be considered in some sense to be the ``regular''
conic representation of $G$. Can any irreducible conic representation
$V$ of $G$, satisfying the requirement that the sections of $\overline{ExtV}$
are strongly proximal, be embedded in $V_{D}$? The answer is negative,
but the minimal factor of $\overline{ExtV}$ (in the category of semi-conic
representations) and the minimal factor of $V$ (in the category of
conic representations) can, in a sense to be made clear, and for this
we do not even require the action on the sections of $\overline{ExtV}$
to be strongly proximal or $V$ to be irreducible. We only require
the induced action on the sections $\overline{ExtV}$ to be $K$-transitive.\\

So let $V$ be a conic representation of $G$ with a $K$-transitive
induced action on the sections of $\overline{ExtV}$. It is proven
in \cite{key-2} that given a compact Hausdorff $G$-space $Y,$ if
the restriction of the action on $Y$ from $G$ to $K$ is transitive
then the natural mapping of the $K$-multipliers sub-group to $H\left(Y\right)$
is a group isomorphism. We know by theorem 30 that the multipliers
on the sections of $\overline{ExtV}$ are exactly all the members
of a coset in $H\left(Y\right)$ and hence include exactly one $K$-multiplier
which we denote by $\sigma$, and we denote by $Y$ a section possessing
it. Hence $\overline{ExtV\left(\sigma\right)}$ is a factor of $\overline{ExtV}$
(incidently this implies that $\overline{ExtV}$ has a factor with
strongly proximal sections) and similarly $V\left(\sigma\right)$
is a factor of $V$ since a multiplier of a section is a $K$-multiplier
if it is a $K$-Multiplier on the closure the extreme points of the
section.\\

Let $\psi:\overline{ExtV}\rightarrow W'$ be any factor of $\overline{ExtV}$.
$K$ acts transitively also on sections of $W'$ and thus $W'$ has
a section $Y'$ whose multiplier is a $K$-multiplier. $\psi^{-1}\left(Y'\right)$
must be a section with the $K$-multiplier of $W$, that is a positive
multiple of $Y$ (by prop. 28). So we can select $Y'$ such that $\psi^{-1}\left(Y'\right)=Y$.
Hence, if $\psi\left(y_{1}\right)=\psi\left(y_{2}\right)$ then $\sigma\left(g,y_{1}\right)=\sigma\left(g,y_{2}\right)$
for any $g\in G$, and that means $\varphi\left(y_{1}\right)=\varphi\left(y_{2}\right)$.
By the universal property of quotient mappings, since $\varphi$ respects
the fibers of $\psi$, it necessarily factors (linearly) through it.
We have thus established the fact that $\overline{ExtV\left(\sigma\right)}$
is the minimal factor of $\overline{ExtV}$. It is possible to show
in an analogous manner that $V\left(\sigma\right)$ is the minimal
factor of $V$. In particular this means $V\left(\sigma\right)$ is
prime.\\

In summary, all conic representations $V$ of $G$ on which the induced
action on sections of $\overline{ExtV}$ is $K$-transitive have minimal
factors (in the sense described in the last paragraph), and the latter
can be embedded in $V_{D}$. However, we do not know if all irreducible
conic representations of $V$ have the mentioned $K$-transitive property.\\

\section{The Case $G=SL_{2}\left(\mathbb{R}\right)$}

As was already mentioned, the action of $SL_{2}\left(\mathbb{R}\right)$
on $P^{1}$ is the universal strongly proximal space of the group,
and the corresponding action on $\Pr\left(P^{1}\right)$ is the universal
irreducible affine representation of the group. It was long-known
that in fact the former is the only non-trivial strongly proximal
space (as mentioned in {[}2{]}), but it was recently shown that the
latter is the only non-trivial irreducible affine representation of
the group (see \cite{key-5}). However, there are non-trivial multipliers
for the action of the group on $P^{1}$. This can be seen since a
necessary condition for a multiplier $\sigma\left(g,x\right)$ to
be trivial is to be dependent on $x$ and $gx$. Taking $m$ the uniform
normalized measure on $P^{1}$ - following \cite{key-2} - we define
$\sigma\left(g,x\right)=\frac{d\left(g^{-1}m\right)}{dm}\left(x\right)$
for all $g\in SL_{2}\left(\mathbb{R}\right)$, $x\in P^{1}$ (the
Radon-Nikodym derivative is positive and continuous), and it can be
checked it is a multiplier by a straight-forward calculation \footnote{For a homogeneous space of a topological group, a measure which is
equivalent (mutual absolute continuity) to all its translations is
called quasi-invariant. If the Radon-Nikodym derivative is continuous
it gives rise to a multiplier in the same fashion.}. For $g=\left(\begin{array}{cc}
a & 1\\
0 & a^{-1}
\end{array}\right)$, $a>>1$, the measure $gm$ obviously has the highest density at
$\overline{\left(1,0\right)}$ and its total mass is $1$, hence $\sigma\left(g,\overline{\left(1,0\right)}\right)=\frac{d\left(g^{-1}m\right)}{dm}\left(\overline{\left(1,0\right)}\right)>1$.
However, for $I=\left(\begin{array}{cc}
1 & 0\\
0 & 1
\end{array}\right)$ $\sigma\left(I,\overline{\left(1,0\right)}\right)=1$, but $g\overline{\left(1,0\right)}=I\overline{\left(1,0\right)}=\overline{\left(1,0\right)}$,
so the necessary condition for $\sigma$ being trivial does not hold.
Actually, for the same reason, the $K$-multipliers - For $K$=$SO\left(2\right)$
- $\sigma^{s}$ for $0\neq s\in\mathbb{R}$ are all non-trivial and
it can be shown the mapping $\mathbb{R}\rightarrow H\left(P^{1}\right)$
that sends $s$ to $\sigma^{s}$ is a group isomorphism (see \cite{key-2}).
So for each isomorphism type of conic representations of $SL_{2}\left(\mathbb{R}\right)$
induced by the strongly proximal action on $P^{1}$ and a multiplier
there exists exactly one $\sigma^{s}$ that induces it.\\

However not every $\sigma^{s}$ induces an irreducible conic representation.
Despite being a multiplier of a strongly proximal action, it was already
noticed in \cite{key-2} that in the case $s=1$, $V\left(\sigma\right)$
is not an irreducible conic representation of the group because $\sigma\left(g,m\right)=\int_{P^{1}}\sigma\left(g,x\right)dm\left(x\right)=1$
is a fixed point. There is a characterization in \cite{key-2} for
when $\sigma^{s}$ gives rise to an irreducible representation $V\left(\sigma^{s}\right)$
in $V_{D}$ ($D$= $SO\left(2\right)$\textbackslash$SL_{2}\left(R\right)$ ).
In a certain way described there, each zonal spherical function on
the group corresponds to exactly one $\sigma^{r}$, and different
spherical functions correspond to different ones (so the correspondence
is one-to-one but not onto). It is proven that the multipliers that
give rise to irreducible representations are precisely $\sigma^{1-r}$
when $\sigma^{r}$ is the corresponding multiplier of some zonal spherical
function. Since $P^{1}$ has no non-trivial factors, all these $V\left(\sigma^{1-r}\right)$
with $r\neq1$ have $P^{1}$ as the section of their $\overline{Ext}$,
and thus are non-isomorphic. Since $SL_{2}\left(\mathbb{R}\right)$
has a continuum of different zonal spherical functions (see \cite{key-8}),
they thus induce a continuum of non-isomorphic irreducible conic representations.
In particular, since all those multipliers but one are non-trivial,
the group has non-degenerate irreducible conic representations and
hence has no universal one. \\

The only other irreducible representation of the group known to us
is the degenerate one-dimensional one. It is interesting to notice
that it is a factor of the degenerate conic representation belonging
to the universal irreducible affine representation and the other irreducible
conic representations mentioned above do not have one-dimensional
factors. In fact, they are all prime and thus have no factors at all.
This can be seen in two independent ways: either by the analysis of
the preceding section (that showed minimality which is stronger) or
by prop. 37. By prop. 37 we also know their sections are simplices,
i.e. for any section $Q$, the barycentric mapping $\beta:\Pr\left(\overline{ExtQ}\right)\rightarrow Q$
is one-to-one (an isomorphism). Up to isomorphism - these are all
the irreducible representations of the group such that the closure
of the extreme points of their sections is isomorphic to $P^{1}$
with the strongly proximal action.\\

It is unknown whether other irreducible conic representations of the
group exist (such representations should necessarily have on the sections
of its $\overline{Ext}$ semi-conic pairs with an action that is not
strongly proximal and a non-trivial multiplier). It is tempting to
guess that the answer to this question is negative since they will
not appear in $V_{D}$, which is a candidate for the ``regular''
conic representation of $SL_{2}\left(\mathbb{R}\right)$. Furthermore,
if such an irreducible $V$ does exist, the group action need not
only be not strongly proximal on sections of $\overline{ExtV}$, but
either $K$ must be non-transitive on them or $\overline{ExtV}$ must
have $\overline{Ext}$ of one of the above $V\left(\sigma^{1-r}\right)$
as its factor. \\

\end{document}